\documentclass[11pt,reqno]{amsart}
\usepackage{amsmath,amssymb,bbm,rotating,times,url}

\newtheorem{theorem}{Theorem}
\newtheorem{lemma}[theorem]{Lemma}
\newtheorem{proposition}[theorem]{Proposition}
\newtheorem{corollary}[theorem]{Corollary}

\newcommand{\e}{\varepsilon}
\newcommand{\charfn}{\mathbbm{1}}
\newcommand{\cube}{{\mathcal C}}
\newcommand{\C}{{\mathbb C}}

\newcommand{\N}{{\mathbb N}}

\newcommand{\R}{{\mathbb R}}
\newcommand{\Rd}{{\R^d}}

\newcommand{\sign}{\operatorname{sign}}
\newcommand{\spt}{\operatorname{spt}}
\newcommand{\widechecknu}{\nu \hspace{-2.5pt} \raisebox{14pt}[0pt][0pt]{\begin{rotate}{180} $\widehat{}$ \end{rotate}} \hspace{3pt}}
\newcommand{\Z}{{\mathbb Z}}
\newcommand{\Zd}{{\Z^d}}

\begin{document}

\title[]{Affine synthesis and coefficient norms for Lebesgue, Hardy and Sobolev spaces}

\author[]{H.-Q. Bui and R. S. Laugesen}
\address{Department of Mathematics, University of Canterbury,
  Christchurch 8020, New Zealand}
\email{Q.Bui\@@math.canterbury.ac.nz}
\address{Department of Mathematics, University of Illinois, Urbana,
IL 61801, U.S.A.} \email{Laugesen\@@uiuc.edu}
\date{\today}



\thanks{Laugesen was supported by N.S.F. Award DMS--0140481 and an
Erskine Visiting Fellowship from the University of Canterbury.}

\begin{abstract}
The affine synthesis operator $Sc=\sum_{j>0} \sum_{k \in \Zd}
c_{j,k} \psi_{j,k}$ is shown to map the mixed-norm sequence space
$\ell^1(\ell^p)$ surjectively onto $L^p(\Rd), 1 \leq p < \infty$,
under mild conditions on the synthesizer $\psi \in L^p(\Rd)$ (say,
having a radially decreasing $L^1$ majorant near infinity) and
assuming $\int_\Rd \psi \, dx = 1$. Here $\psi_{j,k}(x)=|\det
a_j|^{1/p} \psi(a_j x - k)$, for some dilation matrices $a_j$ that
expand.

Hence the standard norm on $f \in L^p(\Rd)$ is equivalent to the
minimal coefficient norm of realizations of $f$ in terms of the
affine system:
\[
\| f \|_p \approx \inf \left\{ \sum_{j>0} ( \sum_{k \in \Zd}
|c_{j,k}|^p )^{1/p} : f=\sum_{j>0} \sum_{k \in \Zd} c_{j,k}
\psi_{j,k} \right\} .
\]

We further show the synthesis operator maps a discrete Hardy space
onto $H^1(\Rd)$, which yields a norm equivalence involving
convolution with a discrete Riesz kernel sequence $\{ z_\ell \}$:
\[
\| f \|_{H^1} \approx \inf \left\{ \sum_{j>0} \sum_{k \in \Zd} (
|c_{j,k}| + |\sum_{\ell \in \Zd} z_\ell c_{j,k-\ell}| ) :
f=\sum_{j>0} \sum_{k \in \Zd} c_{j,k} \psi_{j,k} \right\} .
\]

Coefficient norm equivalences are established also for the Sobolev
spaces $W^{m,p}(\Rd)$, by applying difference operators to the
coefficient sequence $c_{j,k}$.
\end{abstract}

\maketitle

\section{\bf Introduction}
\label{introduction}

This paper studies mapping properties of the \emph{affine
synthesis} operator
\[
c = \{ c_{j,k} \} \mapsto \sum_{j > 0} \sum_{k \in \Zd} c_{j,k}
\psi_{j,k} = Sc
\]
where
\[
\psi_{j,k}(x)=|\det a_j|^{1/p} \psi(a_j x - bk) , \qquad x \in \Rd ,
\]
assuming $\int_\Rd \psi  \, dx = 1$. Affine synthesis arises
naturally in harmonic analysis and approximation theory, as a
discretization of convolution.

We first explain our notation, and then our goals.
\begin{itemize}
\item[\small $\bullet$] The \emph{dimension} $d \in \N$ is fixed
throughout the paper, as is the exponent $p$.

\item[\small $\bullet$]  The \emph{dilation matrices} $a_j$ are
invertible $d \times d$ real matrices that are \emph{expanding}:
\[
\| a_j^{-1} \| \to 0 \qquad \text{as $j \to \infty$.}
\]
(Here $\| \cdot \|$ denotes the norm of a matrix as an operator from
the column vector space $\Rd$ to itself.) For example, one could
take $a_j=2^j I$.

\item[\small $\bullet$] The \emph{translation matrix} $b$ is an
invertible $d \times d$ real matrix, for example the identity
matrix.
\end{itemize}
Our goal is to synthesize surjectively onto the classic function
spaces of analysis, while assuming as little as possible about the
synthesizer $\psi$. This will demonstrate that the ability to
decompose arbitrary functions into linear combinations of the
translates and dilates $\psi_{j,k}$ does not require any special
properties of $\psi$, even though the efficiency of decomposition
naturally does depend on such properties.

\vspace{6pt} \emph{Lebesgue space.} In \cite{bl1} we showed
\[
S : \{ \text{finite sequences} \} \to L^p \quad \text{has dense
range}
\]
for $1 \leq p < \infty$, assuming only that $\psi$ has periodization
locally in $L^p$ (meaning $\sum_{k \in \Zd} |\psi(x-bk)| \in
L^p_{loc}$) and has nonzero integral ($\int_\Rd \psi \, dx \neq 0$).
This density of the range of $S$ means that the small-scale affine
system $\{ \psi_{j,k} : j>0, k \in \Zd \}$ \emph{spans} $L^p$.
Notice Strang--Fix conditions are not imposed: the translates $\{
\psi(\cdot -bk) \}_{k \in \Zd}$ need not form a partition of unity.
And recall from \cite{bl1} that the periodization assumption on
$\psi$ is easily satisfied, holding for example if $\psi \in L^p$
has a radially decreasing $L^1$ majorant near infinity.

Two natural questions arise from this $L^p$ density result: what is
the right \emph{domain} for $S$? and does $S$ map this domain
\emph{onto} $L^p$? We answer these questions in
Section~\ref{lebesgueresults}, where Theorems~\ref{lebesguerep1} and
\ref{lebesguerep2} show that
\[
S : \ell^1(\ell^p) \to L^p \quad \text{is linear, bounded and onto.}
\]
The domain $\ell^1(\ell^p)$ is the mixed-norm space of coefficients
satisfying
\[
\sum_{j>0} ( \sum_{k \in \Zd} |c_{j,k}|^p )^{1/p} < \infty .
\]

From the surjectivity of $S$ onto $L^p$ one deduces a coefficient
norm equivalence of the form
\[
\| f \|_p \approx \inf \{ \| c \|_{\ell^1(\ell^p)} : f=Sc \} .
\]
(Bruna \cite[Theorem~4]{B05} earlier proved the case $p=1$.) The
constants for this norm equivalence are evaluated in
Corollary~\ref{lebesguenorm1}, and Corollary~\ref{lebesguenorm2}
shows that in fact equality holds under suitable conditions.

Corollary~\ref{lebesguenorm4} proves a \emph{localized} norm
equivalence, on a domain $\Omega \subset \Rd$.

\vspace{6pt} \emph{Hardy space.} Section~\ref{hardyresults}
considers the Hardy space $H^1 = H^1(\Rd) = \{ f \in L^1 : f
* (x/|x|^{d+1}) \in L^1 \}$. Again we ask what the domain of the synthesis operator
should be, and whether $S$ maps this domain onto $H^1$.
Theorems~\ref{hardyrep1} and \ref{hardyrep2} answer the question by
showing
\[
S : \ell^1(h^1) \to H^1 \quad \text{is linear, bounded and onto,}
\]
provided $\psi \in L^1$ has nonzero integral and is somewhat
``nice'' (being for example an $L^p$ function with compact support
for some $p>1$, or a Schwartz function). The discrete Hardy space
$h^1$ here is defined in Section~\ref{hardyresults} by convolution
of sequences against a ``singular'' kernel in the $k$-variable
(which is the discrete analogue of convolution against $x/|x|^{d+1}$
in the definition of the continuous Hardy space $H^1$).

Notice the synthesizer $\psi$ cannot belong to $H^1$, because it has
nonzero integral, and so to ensure the linear combination $Sc$
belongs to $H^1$ we must invoke cancellation properties of the
coefficient space $h^1$. This approach to synthesis in $H^1$ seems
natural to us because the analysis operator takes $H^1$ to $h^1$
(see below), and one would like to reconstruct $H^1$ using only the
data obtained from analysis.

Corollary~\ref{hardynorm} deduces a coefficient norm equivalence
of the form
\[
\| f \|_{H^1} \approx \inf \{ \| c \|_{\ell^1(h^1)} : f=Sc \} .
\]

\vspace{6pt} \emph{Sobolev space.} For synthesis into the Sobolev
space $W^{m,p}(\Rd)$, we establish boundedness and surjectivity in
Theorems~\ref{sobolevrep1} and \ref{sobolevrep2} of
Section~\ref{sobolevresults}. The sequence space from which we
synthesize involves difference operators with respect to the
$k$-index, which act as discrete analogues of differentiation.

\vspace{6pt} \emph{Discussion.} Our results on surjectivity of the
synthesis operator seem to be qualitatively new --- they rely on a
method of ``scale-averaged convergence'' that we developed only
recently for $L^p$ in \cite{bl1}. Boundedness of the synthesis
operator has of course been studied before.

We also consider boundedness of the affine \emph{analysis operator}
at each scale $j$, in other words, boundedness of $f \mapsto \{
\langle f , \phi_{j,k} \rangle \}_{k \in \Zd}$ from $L^p$ to
$\ell^p$, from $H^1$ to $h^1$ and from $W^{m,p}$ to $w^{m,p}$. See
Propositions~\ref{analysisleb}, \ref{analysishardy} and
\ref{analysissob} respectively. We show the full analysis operator
(over all scales) maps isomorphically onto its range in
Corollaries~\ref{lebesguenorm3}, \ref{hardynormb} and
\ref{sobolevnorm3}, thus giving coefficient norms in terms of the
analysis operator.

If one works at a fixed scale $j$, rather than considering all
scales $j>0$ as we do in this paper, then synthesis yields a shift
invariant subspace of $L^p$. Aldroubi, Sun and Tang's $p$-frame work
for such shift invariant spaces \cite{AST01} is described at the end
of Section~\ref{lebesgueresults}.

Topics we do not pursue in this paper include Gabor systems
(modulations and translations) and wavepacket decompositions
(modulations, translations and dilations). For some recent work in
these areas one can consult \cite{CKS06,FF06,FES98,Gro01,LWW04}.

The paper is structured as follows. Section~\ref{definition}
establishes notation and definitions.
Sections~\ref{lebesgueresults}--\ref{sobolevresults} present our
synthesis results on Lebesgue, Hardy and Sobolev spaces. The proofs
are in Sections~\ref{lebesgue_proofs}--\ref{sobolev_proofs}.
Appendices treat discrete Hardy spaces, and Banach frames.

Parallel results on Triebel--Lizorkin spaces are discussed in
Section~\ref{literature}. Open problems when $\psi$ has zero
integral ($\int_\Rd \psi \, dx = 0$) are treated in
Section~\ref{openproblems}, including Meyer's Mexican hat spanning
problem for $L^p$ and its counterpart for the Hardy space. We hope
this paper helps contribute towards an eventual resolution of these
fascinating open problems.

\section{\bf Further definitions and notation} \label{definition}

1. We use doubly-indexed sequences $c=\{ c_{j,k} \}_{j>0,k \in
\Zd}$ of complex numbers, with the norm
\[
\| c \|_{\ell^1(\ell^p)} := \sum_{j>0} \left( \sum_{k \in \Zd}
|c_{j,k}|^p \right)^{\! \! \! 1/p}
\]
when $1 \leq p < \infty$. When $p=\infty$, define $\| c
\|_{\ell^1(\ell^\infty)} := \sum_{j>0} \sup_{k \in \Zd}
|c_{j,k}|$. Then
\[
\ell^1(\ell^p) := \left\{ c : \| c \|_{\ell^1(\ell^p)} < \infty
\right\}
\]
is a Banach space.

2. Write $L^p=L^p(\Rd)$ for the class of complex valued functions
with finite $L^p$-norm. Given $\psi \in L^p$ and $\phi \in L^q$,
with
\[
\frac{1}{p}+\frac{1}{q}=1
\]
by notational convention, we define rescalings
\[
\psi_{j,k}(x) = |\det a_j|^{1/p} \psi(a_j x-bk) , \quad
\phi_{j,k}(x) = |\det a_j|^{1/q} \phi(a_j x-bk) .
\]
These rescalings preserve the $L^p$-norm $\| \psi_{j,k} \|_p = \|
\psi \|_p$ and the $L^q$-norm $\| \phi_{j,k} \|_q = \| \phi \|_q$,
respectively. \emph{Alert:} the notation $\psi_{j,k}$ conceals its
dependence on $p$.

3. The \emph{synthesis operator} is
\begin{align*}
Sc = S_{\psi,b} c = \sum_{j>0} \sum_{k \in \Zd} c_{j,k} \psi_{j,k} .
\end{align*}
Our theorems will specify acceptable domains for this operator, and
will explain the sense in which the sums over $j$ and $k$ converge.
Occasionally we will synthesize at a fixed scale $j$ by writing
\[
S_j s = \sum_{k \in \Zd} s_k \psi_{j,k} ,
\]
for sequences $s=\{ s_k \}_{k \in \Zd}$.

4. The \emph{analysis operator at scale $j$} is
\[
T_j f = T_{j,\phi} f = \{ |\det b| \langle f , \phi_{j,k} \rangle
\}_{k \in \Zd} .
\]
That is, $T_j$ maps a function $f$ to its sequence of sampled
$\phi$-averages at scale $j$. The full \emph{analysis operator}
simply combines these sequences as
\[
Tf = T_\phi f = \{ |\det b| \langle f , \phi_{j,k} \rangle \}_{j>0,
k \in \Zd} .
\]

Our analysis and synthesis operators depend implicitly on the
exponent $p$, through the definitions of $\phi_{j,k}$ and
$\psi_{j,k}$.

5. The \emph{periodization} of a function $f$ is
\[
Pf(x) = |\det b| \sum_{k \in \Zd} f(x-bk) \qquad \text{for $x \in
\Rd$.}
\]

6. Write $\cube=[0,1)^d$ for the unit cube in $\Rd$, and
$\cube_0=(-1/2,1/2)^d$ for the centered open unit cube. We regard
$\cube$ as consisting of column vectors and $\cube_0$ as
consisting of row vectors, as the context will always make clear.

\section{\bf $L^p$ results} \label{lebesgueresults}

First we obtain boundedness of the synthesis operator, when the
periodization of the synthesizer belongs locally to $L^p$. This was
already observed by Aldroubi, Sun and Tang \cite[formula
(2.3)]{AST01}, but we give a proof in
Section~\ref{lebesguerep1_proof} anyway, to keep the paper
self-contained.
\begin{theorem}[Synthesis into $L^p$] \label{lebesguerep1}
Assume $1 \leq p \leq \infty$ and $\psi \in L^p$ with $P|\psi| \in
L^p_{loc}$.

Then $S : \ell^1(\ell^p) \to L^p$ is bounded. More precisely, if $c
\in \ell^1(\ell^p)$ then the series
\begin{equation} \label{frep}
Sc = \sum_{j>0} \sum_{k \in \Zd} c_{j,k} \psi_{j,k}
\end{equation}
converges in $L^p$ in the sense that
\begin{equation} \label{lebesgueconverge}
\begin{minipage}{.9\linewidth}
the sum over $k$ in \eqref{frep} converges pointwise absolutely
a.e.\ to a function in $L^p$, and the sum over $j$ converges
absolutely in $L^p$,
\end{minipage}
\end{equation}
and furthermore
\begin{equation} \label{analogue1}
\| Sc \|_p \leq |\det b|^{-1} \| P|\psi| \|_{L^p(b\cube)} \| c
\|_{\ell^1(\ell^p)} .
\end{equation}
\end{theorem}
After proving the theorem in Section~\ref{lebesguerep1_proof}, we
also give examples to show $\ell^1(\ell^p)$ is the ``correct''
domain for $S$, when the synthesizer $\psi$ has nonzero integral.
When $p=2$ and $\psi$ has zero integral, we point out that synthesis
can be bounded on the larger domain $\ell^2(\ell^2)$, for wavelet
and affine frame generators.

Remember that $Sc$ in \eqref{frep} depends implicitly on the value
of $p$, through the renormalization factor $|\det a_j|^{1/p}$ in the
definition of $\psi_{j,k}$. This dependence would be problematic if
we synthesized into more than one $L^p$-space at a time, but we will
not.

The hypothesis that the periodization of $|\psi|$ be locally in
$L^p$ is rather weak, and is easily verified in many cases. For
example when $p=1$ it holds for all $\psi \in L^1$. And for $p>1$ it
holds when $\psi \in L^p$ has compact support or when $\psi$ has a
bounded, radially decreasing $L^1$-majorant, or when $\psi$ equals a
sum of such functions. On the other hand, $P|\psi| \in L^p_{loc}$
can hold even when $\psi$ does not decay at infinity. See
\cite[\S3.1]{bl1} for all these observations. Thus
Theorem~\ref{lebesguerep1} improves somewhat on earlier boundedness
results (which go back as far as \cite{G69,SF73,S93}) because $\psi$
need have neither compact support nor decay at infinity.

In the case $p=2$, the periodization hypothesis on $\psi$ can be
weakened to just $\psi \in L^2$ with $P(|\widehat{\psi}|^2) \in
L^\infty$. Specifically, one has bounded synthesis with $\| Sc \|_2
\leq \| P(|\widehat{\psi}|^2) \|_\infty^{1/2} \| c
\|_{\ell^1(\ell^2)}$ where this last periodization is taken with
respect to the lattice $\Zd b^{-1}$.
This estimate is proved in \cite[Theorem~7.2.3]{Ch03}, by showing
that $P(|\widehat{\psi}|^2) \in L^\infty$ is exactly the condition
for the integer translates of $\psi$ to form a Bessel sequence (that
is, to satisfy an upper frame bound).

For all $p$ one might ask whether the periodization assumption on
$|\psi|$ in Theorem~\ref{lebesguerep1} can be weakened to just $\psi
\in L^1 \cap L^p$. We do not know.

\vspace{6pt} Next we show the synthesis operator is surjective, when
$1 \leq p < \infty$. In other words, we show every $f \in L^p$ can
be expressed by a series of the form \eqref{frep}.
\begin{theorem}[Synthesis onto $L^p$] \label{lebesguerep2}
Assume $1 \leq p < \infty$ and $\psi \in L^p$ with $P|\psi| \in
L^p_{loc}$ and $\int_\Rd \psi \, dx = 1$.

Then $S : \ell^1(\ell^p) \to L^p$ is open, and surjective. Indeed,
if $f \in L^p$ and $\e>0$ then a sequence $c \in \ell^1(\ell^p)$
exists such that $Sc=f$ with convergence as in
\eqref{lebesgueconverge}, and such that
\[
\| c \|_{\ell^1(\ell^p)} \leq |\det b|^{1/q} \| f \|_p + \e.
\]
\end{theorem}
Section~\ref{lebesguerep2_proof} has the proof. The integral of
$\psi$ is well defined, in the statement of
Theorem~\ref{lebesguerep2}, because the assumption $P|\psi| \in
L^p_{loc}$ implies $P|\psi| \in L^1_{loc}$ and hence $\psi \in L^1$.

We are not aware of any prior general work on surjectivity of the
synthesis operator, when $\psi$ has nonzero integral. The closest
work seems to be Filippov and Oswald's construction in
\cite[Theorems~1 and 3]{FO95}, \cite{F98}, of ``representation
systems'' by which every $f \in L^p$ can be written as a convergent
series $Sc=f$. This looks like surjectivity, but the drawback is
that their result yields no control over the size of coefficients in
the sequence $c$, and thus it is unclear what the domain of $S$
really is. Further, Filippov and Oswald work only with isotropic
dilation matrices.

We also mention the density results coming from Strang--Fix theory
(discussed in \cite[\S3]{bl1}, although note
Theorem~\ref{lebesguerep2} holds without needing the Strang--Fix
hypotheses). When $\psi$ has zero integral, the wavelet theory
\cite{Ch03,D92,HW96,M92} and related phi-transform theory (described
in Section~\ref{literature}) provide sufficient conditions for
obtaining frames, orthonormal bases, and unconditional bases,
provided the large scales $j \leq 0$ are included in the synthesis.
These conditions all ensure surjectivity of $S$ on suitable domains.
But the zero-integral case also raises intriguing open problems,
discussed in Section~\ref{openproblems}.

For the special case $p=2$, we remarked above that the condition
$P(|\widehat{\psi}|^2) \in L^\infty$ implies bounded synthesis. We
suspect it also implies surjectivity onto $L^2$, provided
$|\widehat{\psi}|$ is continuous near the origin and
$\widehat{\psi}(0) = 1$. These conditions certainly guarantee the
$\psi_{j,k}$ span $L^2$, by a result of Daubechies
\cite[Proposition~5.3.2]{D92}, and thus $S$ has dense range. One
would like to improve this to full range.

\vspace{6pt} \noindent \emph{Remark on non-injectivity.} The
synthesis operator is not injective, and indeed has a very large
kernel. For example, we could discard the dilation $a_1$ (in other
words, discard all terms with $j=1$ in the sum defining $Sc$) and
still show $S$ maps onto $L^p$, by applying
Theorem~\ref{lebesguerep2} with the remaining dilations $\{
a_2,a_3,\ldots \}$.

\vspace{6pt} Equivalence of the $L^p$ and $\ell^1(\ell^p)$ norms
follows immediately from Theorems~\ref{lebesguerep1} and
\ref{lebesguerep2}:
\begin{corollary}[Synthesis norm for $L^p$] \label{lebesguenorm1}
Assume $1 \leq p < \infty$ and $\psi \in L^p$ with $P|\psi| \in
L^p_{loc}$ and $\int_\Rd \psi \, dx = 1$. Then
\[
\| f \|_p \approx \inf \left\{ \| c \|_{\ell^1(\ell^p)} : \text{$f =
Sc$ \, as in \eqref{lebesgueconverge}} \right\}
\]
for all $f \in L^p$. Explicitly,
\begin{align*}
|\det b| \| P|\psi| \|_{L^p(b\cube)}^{-1} \| f \|_p & \leq \inf
\left\{ \| c \|_{\ell^1(\ell^p)} : \text{$f = Sc$ \, as in
\eqref{lebesgueconverge}} \right\} \\
& \leq |\det b|^{1/q} \| f \|_p .
\end{align*}
\end{corollary}
For $p=1$, the corollary was proved by Bruna \cite[Theorem~4]{B05}.
His duality methods apply without needing our assumption that the
translations be restricted to a lattice. Recall when $p=1$ that the
periodization condition $P|\psi| \in L^p_{loc}$ is superfluous,
holding automatically for all $\psi \in L^1$.

\vspace{6pt} \noindent \emph{Remark on norm equivalence.} As soon as
$S$ maps \emph{onto} $L^p$, we know the map $\widetilde{S} :
\ell^1(\ell^p)/\ker S \to L^p$ is a bounded linear bijection. (Here
we use the canonical norm on the quotient space: $\inf_{c^\prime \in
\ker S} \| c + c^\prime \|_{\ell^1(\ell^p)}$.) Then the inverse map
$\widetilde{S}^{-1}$ is bounded by the closed graph theorem, giving
equivalence of the $\ell^1(\ell^p)$ and $L^p$ norms like in
Corollary~\ref{lebesguenorm1}. But Corollary~\ref{lebesguenorm1}
goes further, for it provides an \emph{explicit} upper bound for the
norm equivalence, based on explicitly estimating the norm of
$\widetilde{S}^{-1}$ in Theorem~\ref{lebesguerep2}.

Corollary~\ref{lebesguenorm2} goes further still, giving actual
equality of the norms when $\psi$ is nonnegative.
\begin{corollary}[Synthesis norm equality] \label{lebesguenorm2}
Assume $1 \leq p < \infty$ and that $\psi \in L^p$ is nonnegative.
When $p=1$ assume $\int_\Rd \psi \, dx = 1$, and when $1<p<\infty$
assume $P\psi \equiv 1$ (which implies $\int_\Rd \psi \, dx = 1$).

Then for all $f \in L^p$,
\[
\| f \|_p = |\det b|^{-1/q} \inf \left\{ \| c \|_{\ell^1(\ell^p)} :
\text{$f=Sc$ \, as in \eqref{lebesgueconverge}} \right\} .
\]
\end{corollary}
The constant periodization condition $P\psi \equiv 1$ in this
corollary says that the collection $\{ |\det b| \psi(x - bk) : k \in
\Zd \}$ of translates of $\psi$ is a partition of unity. Examples of
such $\psi$ (when $b=I$) include the indicator function
$\charfn_\cube$ and convolutions of this indicator function with any
nonnegative function having integral $1$.

We next ``localize'' Corollary~\ref{lebesguenorm1} to an open set
$\Omega \subset \Rd$. To that end, we say a sequence $c = \{ c_{j,k}
\}_{j>0, k\in \Zd}$ is \emph{adapted to $\Omega$ and $\psi$} if
$\spt(\psi_{j,k}) \subset \Omega$ whenever $c_{j,k} \neq 0$, or in
other words if $c_{j,k}=0$ whenever $\spt(\psi_{j,k}) \cap \Omega^c
\neq \emptyset$. The point of this definition is to ensure $Sc=0$ on
the complement of $\Omega$.
\begin{corollary}[Synthesis norm for $L^p(\Omega)$] \label{lebesguenorm4}
Assume $\Omega \subset \Rd$ is open and nonempty, take $1 \leq p <
\infty$, and suppose $\psi \in L^p(\Rd)$ is compactly supported
with $\int_\Rd \psi \, dx = 1$. Then for all $f \in L^p(\Omega)$,
\[
\| f \|_{L^p(\Omega)} \approx \inf \left\{ \| c \|_{\ell^1(\ell^p)}
: \text{$f = Sc$ \, as in \eqref{lebesgueconverge}, and $c$ is
adapted to $\Omega$ and $\psi$} \right\} .
\]
\end{corollary}
The constants in this norm equivalence are the same as in
Corollary~\ref{lebesguenorm1}, and so they depend on $\psi$ and $b$
but are independent of $\Omega$. The corollary is proved in
Section~\ref{lebesguenorm4_proof}.

We turn now to the analysis operator, which also yields a
coefficient norm.
\begin{corollary}[Analysis norm for $L^p$]\label{lebesguenorm3}
Assume $1\leq p <\infty$, and take an analyzer $\phi\in L^q$ with
$P|\phi|\in L^\infty$ and $\int_{\Rd} \phi\, dx = 1$.

Then for all $f\in L^p$,
\[
\|f\|_p \approx \| Tf \|_{\ell^\infty(\ell^p)} = \sup_{j>0} \left(
\sum_{k\in\Zd} | \langle f , \phi_{j,k} \rangle|^p \right)^{\! \!
1/p} |\det b|.
\]
\end{corollary}
In other words, the analysis operator is linear, bounded and
injective from $L^p$ onto its range in the mixed norm sequence space
$\ell^\infty(\ell^p)$. To explain the appearance of
$\ell^\infty(\ell^p)$ in the corollary, note the analysis and
synthesis operators are adjoints, with $T_\phi : L^p \to
\ell^\infty(\ell^p)$ being the adjoint of $S_\phi : \ell^1(\ell^q)
\to L^q$, at least when $1 < p \leq \infty$. Thus the injectivity of
analysis in Corollary~\ref{lebesguenorm3} is equivalent to the
surjectivity of synthesis in Theorem~\ref{lebesguerep2}. But we
prove the corollary in Section~\ref{lebesguenorm3_proof} anyway, for
the sake of concreteness and to handle $p=1$.

\vspace{6pt} We close the section by describing Aldroubi, Sun and
Tang's work on $p$-frames \cite{AST01}, which is close in subject
matter to this paper but has little direct overlap. They study the
shift invariant range space $V_j = S_j(\ell^p)$ of the synthesis
operator at a \emph{single} scale $j$ (whereas our results combine
all scales $j>0$). Roughly, they showed that the reconstruction
formula $S_j T_j = \text{identity}$ on $V_j$ (which says functions
in $V_j$ can be synthesized from their sampled average values at
scale $j$) holds if and only if $V_j$ is closed in $L^p$, if and
only if $T_j$ is injective on $V_j$ when $\phi=\psi$, if and only if
$\widehat{\psi}$ satisfies a certain ``bracket product'' condition.
This is all carried out in the multiply generated case, with
synthesizers $\psi_1,\ldots,\psi_r$. They observe that the range
$V_j$ need not be closed, for example when
$\psi=\charfn_{[0,1)}-\charfn_{[1,2)}$ in one dimension with
$1<p<\infty$ \cite[page~7]{AST01}.

\section{\bf Hardy space results} \label{hardyresults}

Our Hardy space results assume the dilation matrices are isotropic
and expanding. To be precise, all the results in this section assume
that
\[
a_j = \alpha_j I
\]
for some ``dilation'' sequence $\alpha = \{ \alpha_j \}_{j>0}$ of
nonzero real numbers with $|\alpha_j| \to \infty$ as $j \to \infty$.

We will recall the Hardy space $H^1$, and then construct a
discrete Hardy space on which the synthesis operator can act. Then
we state our synthesis results.

\subsection*{Hardy space $H^1$.}
Define the Fourier transform with $2\pi$ in the exponent:
\[
\widehat{f}(\xi) = \int_\Rd f(x) e^{-2\pi i \xi x} \, dx,
\]
for row vectors $\xi \in \Rd$. Write $C_d = \Gamma((d+1)/2)
\pi^{-(d+1)/2}$ and
\[
Z(x)=C_d
\begin{cases}
x/|x|^{d+1}, & x \neq 0 ,\\
0 ,  & x = 0 ,
\end{cases}
\]
for the Riesz kernel, so that the Riesz transform of $f \in L^1$
is
\[
Rf(x) = (f*Z)(x) = \text{p.v.} \int_\Rd f(x-y) Z(y) \, dy .
\]
Then $Rf$ is finite a.e., and is a measurable vector-valued
function of $x \in \Rd$. Notice
\[
\widehat{Rf}(\xi) = -i \frac{\xi}{|\xi|} \widehat{f}(\xi).
\]
Recall the \emph{Hardy space} is
\[
H^1 = H^1(\Rd) = \{ f \in L^1 : Rf \in L^1 \}, \quad \text{with the
norm $\| f \|_{H^1} := \| f \|_1 + \| Rf \|_1$.}
\]

Functions in the Hardy space have vanishing integral: if $f \in
H^1$ then $Rf \in L^1$ and so $\widehat{Rf}$ is continuous, which
implies
\begin{equation} \label{intzero}
\widehat{f}(0)=\int_\Rd f(x)\, dx = 0 \quad \text{and} \quad
\widehat{Rf}(0)=\int_\Rd Rf(x)\, dx = 0 .
\end{equation}

The Riesz transform commutes with dilations and translations,
meaning: $R(f(\alpha x-x_0))=\sign(\alpha) (Rf)(\alpha x-x_0)$
when $\alpha \in \R \setminus \{ 0 \}, x_0 \in \Rd$. But dilation
invariance fails when $\alpha$ is an arbitrary matrix, which is
why we restrict to isotropic dilations in this section.

All these facts about Riesz transforms and $H^1$ can be found in
\cite{st1,st2}.

\subsection*{Discrete Hardy space $h^1$.} Take a smooth, compactly
supported cut-off function $\nu$ supported in the centered unit cube
$\cube_0$, with $\nu \equiv 1$ near the origin. Then define a
``discrete Riesz kernel'' sequence $z = \{ z_k \}_{k \in \Zd} \in
\ell^2$ by specifying its Fourier series:
\begin{equation} \label{zdefn}
\zeta(\xi) = \sum_{k \in \Zd} z_k e^{-2\pi i \xi k} := -i \frac{\xi
b^{-1}}{|\xi b^{-1}|} \nu(\xi) , \qquad \xi \in \cube_0 ,
\end{equation}
where $\xi k$ denotes the dot product (recall $\xi$ is a row and $k$
is a column vector) and where for later convenience we use $-\xi$
rather than $+\xi$ in the exponent of the Fourier series.

The sequence $z$ is vector-valued (since $z_k \in \C^d$), and
belongs to $\ell^2$ because $(-i \xi b^{-1}/|\xi b^{-1}|) \nu(\xi)$
is bounded and hence belongs to $L^2(\cube_0)$. Thus the series for
the periodic function $\zeta$ converges unconditionally in
$L^2(\cube_0)$.
When $b=I$, observe from \eqref{zdefn} that $\zeta$ is simply a
cut-off version of the Fourier transform of the Riesz kernel.

Naturally $s * z \in \ell^2$ whenever $s \in \ell^1$. We define a
``discrete Hardy space'' by requiring that $s * z$ belong to the
smaller space $\ell^1$:
\[
h^1 = \left\{ s \in \ell^1 : s * z \in \ell^1 \right\} ,
\]
with a norm
\[
\| s \|_{h^1} := \| s \|_{\ell^1} + \| s * z \|_{\ell^1}
\]
that makes $h^1$ a Banach space.

Appendix~\ref{hardyapp} investigates some properties of $h^1$,
including its independence from the cut-off function $\nu$, and its
relation to the atomic sequence space $H^1(\Zd)$ studied by several
authors. The appendix also points out that sequences in $h^1$ have
vanishing mean: $\sum_{k \in \Zd} s_k = 0$.

\vspace{6pt} The mixed-norm Banach space we will need is
\begin{align*}
\ell^1(h^1) & = \left\{ c \in \ell^1(\ell^1) : c * z \in
\ell^1(\ell^1) \right\} \\
& = \{ c \in \ell^1(\ell^1) : \| c \|_{\ell^1(h^1)} < \infty \} ,
\end{align*}
where the convolution is taken with respect to the $k$-index and
the norm is
\begin{align*}
\| c \|_{\ell^1(h^1)} & := \| c \|_{\ell^1(\ell^1)} + \| c * z
\|_{\ell^1(\ell^1)} \\
& = \sum_{j>0} \| c_j \|_{h^1} = \sum_{j>0} \sum_{k \in \Zd} \left(
|c_{j,k}| + |(c
* z)_{j,k}| \right)
\end{align*}
with the notation $c_j=\{ c_{j,k} \}_{k \in \Zd}$.

\subsection*{Properties of the synthesis operator}
We will prove boundedness and surjectivity of the synthesis
operator. Note
\[
\psi_{j,k}(x) = |\det a_j| \psi(a_j x - bk)
\]
in the next two theorems, because we implicitly take $p=1$ when
working with $\psi \in L^1$.

First we show boundedness.
\begin{theorem}[Synthesis into $H^1$] \label{hardyrep1}
Assume $\psi \in L^1$ and
\begin{equation} \label{h1bound}
\sup_{|y| \leq 1} \| \psi - \psi(\cdot- y) \|_{H^1} < \infty .
\end{equation}

Then $S : \ell^1(h^1) \to H^1$ is bounded. More precisely, if $c \in
\ell^1(h^1)$ then the series $Sc = \sum_{j>0} \sum_{k \in \Zd}
c_{j,k} \psi_{j,k}$ converges in $H^1$ in the sense that
\begin{equation} \label{hardyconverge}
\begin{minipage}{.9\linewidth}
the sum over $k$ converges absolutely in $L^1$ to a function
belonging to $H^1$ and the sum over $j$ converges absolutely in
$H^1$,
\end{minipage}
\end{equation}
and furthermore $\| Sc \|_{H^1} \leq C \| c \|_{\ell^1(h^1)}$ for
some constant $C=C(\psi,b)$.
\end{theorem}
See Section~\ref{hardyrep1_proof} for the proof. To understand why
we so carefully describe the convergence of $Sc$ in $H^1$, in
\eqref{hardyconverge}, just recall that $\psi_{j,k} \notin H^1$ when
$\int_\Rd \psi \, dx \neq 0$.

Many functions $\psi$ satisfy the finite supremum assumption in
\eqref{h1bound}, for example if $\psi \in L^p$ for some $p>1$ and
$\psi$ has compact support, or if $\psi$ is a Schwartz function;
cf.\ \cite[\S3.3-3.4]{bl3}. Incidentally, the supremum in
\eqref{h1bound} can equivalently be taken over any other ball of
$y$-values, by \cite[Lemma~6]{bl3}.

Theorem~\ref{hardyrep1} was proved for synthesizers $\psi \in L^2$
having compact support by S. Boza and M. J. Carro in
\cite[Proposition~3.11]{BC98}, \cite[Theorem~3.1]{BC02}. Their
methods are very different from ours, involving a maximal
characterization of $H^p(\Zd), 0 < p \leq 1$.
Theorem~\ref{hardyrep1} was also proved earlier in one dimension for
synthesizers $\psi \in L^1(\R)$ having compact support and locally
integrable Hilbert transform $R \psi \in L^1_{loc}(\R)$, by Q. Sun
\cite[Theorem 13]{S93}. Note that Sun's assumptions on $\psi$ imply
the condition (\ref{h1bound}). Interestingly, Sun proved a converse
theorem for compactly supported $\psi$, saying that bounded
synthesis implies $R\psi$ must be locally integrable.

Now we show surjectivity of synthesis from $\ell^1(h^1)$ to $H^1$, a
result that seems to be qualitatively new.
\begin{theorem}[Synthesis onto $H^1$] \label{hardyrep2}
Assume $\psi \in L^1$ with $\int_\Rd \psi \, dx = 1$ and
\begin{equation} \label{h1limit}
\| \psi - \psi(\cdot- y) \|_{H^1} \to 0 \qquad \text{as $y \to 0$.}
\end{equation}

Then $S : \ell^1(h^1) \to H^1$ is open, and surjective. Indeed if $f
\in H^1$ and $\e>0$ then a sequence $c \in \ell^1(h^1)$ exists such
that $Sc=f$ (with convergence as in \eqref{hardyconverge}) and
\[
\| c \|_{\ell^1(h^1)} \leq C \| f \|_{H^1} + \e ,
\]
for some constant $C=C(b)$.
\end{theorem}
See Section~\ref{hardyrep2_proof} for the proof. The constant $C$
can be evaluated, if desired.

Assumption \eqref{h1limit} holds, for example, if $\psi \in L^p$ for
some $p>1$ and $\psi$ has compact support, or if $\psi$ is a
Schwartz function (cf.\ \cite[\S3.3]{bl3}).

By combining the last two theorems, we obtain a norm for $H^1$ in
terms just of coefficients in affine expansions.
\begin{corollary}[Synthesis norm for $H^1$] \label{hardynorm}
Assume $\psi \in L^1$ with $\int_\Rd \psi \, dx = 1$, and suppose
\eqref{h1limit} holds. Then for all $f \in H^1$,
\[
\| f \|_{H^1} \approx \inf \left\{ \| c \|_{\ell^1(h^1)} : \text{$c
\in \ell^1(h^1)$ and $f = Sc$ as in \eqref{hardyconverge}} \right\}
.
\]
\end{corollary}
The proof is easy: assumption \eqref{h1limit} implies
\eqref{h1bound} by \cite[\S3.4]{bl3}, and so
Theorems~\ref{hardyrep1} and \ref{hardyrep2} both apply here, giving
the corollary.

The analysis operator also provides a coefficient norm for $H^1$:
\begin{corollary}[Analysis norm for $H^1$]\label{hardynormb}
Take an analyzer $\phi \in L^1$ with $\int_{\Rd} \phi\,dx = 1$, and
suppose $\phi$ is a Schwartz function with $\widehat{\phi}$
supported in $\cube_0 b^{-1}$.

Then for all $f\in H^1$,
\[
\| f \|_{H^1} \approx \| Tf \|_{\ell^\infty(h^1)} = \sup_{j > 0} \|
\{ \langle f,\phi_{j,k} \rangle \} \|_{h^1} \, |\det b| .
\]
\end{corollary}
The corollary says the analysis operator $T$ is linear, bounded and
injective from $H^1$ onto its range in $\ell^\infty(h^1)$. See
Section~\ref{hardynormb_proof} for the proof.

\section{\bf Sobolev space results} \label{sobolevresults}

Write $W^{m,p}=W^{m,p}(\Rd)$ for the class of Sobolev functions
with $m$ derivatives in $L^p$, normed by
\[
\| f \|_{W^{m,p}} = \sum_{|\rho| \leq m} \| D^\rho f \|_p .
\]
Here $\rho=(\rho_1,\ldots,\rho_d)$ is a multiindex of order
$|\rho|=\rho_1+\cdots+\rho_d$.

We continue to assume the dilation matrices are isotropic and
expanding, with $a_j = \alpha_j I$, like in the previous section.

We first construct a sequence space on which the synthesis
operator will act, and then construct a class of synthesizers,
before stating our Sobolev synthesis results.

\vspace{6pt} \emph{Discrete Sobolev space.} Define difference
operators on sequences $s=\{ s_k \}_{k \in \Zd}$ by
\[
\Delta_t s = \{ s_k - s_{k-e_t} \}_{k \in \Zd} ,
\]
for $t=1,\ldots,d$, where $e_t$ is the unit vector in the $t$-th
coordinate direction. Define higher difference operators by
\[
\Delta^\rho s = \Delta_1^{\rho_1} \cdots \Delta_d^{\rho_d} s .
\]
Then we can define a discrete Sobolev space by
\[
w^{m,p} = \left\{ s \in \ell^p : \text{$\Delta^\rho s \in \ell^p$
for each multiindex $\rho$ of order $|\rho| \leq m$} \right\} ,
\]
with norm
\[
\| s \|_{w^{m,p}} = \sum_{|\rho| \leq m} \| \Delta^\rho s
\|_{\ell^p} .
\]
This space is just $\ell^p$ with a new norm, but we proceed to
weight the norm by appropriate powers of the dilation sequence, as
follows.

Recall sequences are multiplied term-by-term, meaning $(c
\tilde{c})_{j,k} = c_{j,k} \tilde{c}_{j,k}$ and so on. In particular
we can define powers $\alpha^r = \{ \alpha_j^r \}_{j>0}$ of the
dilation sequence, whenever $r$ is a nonnegative integer. With these
conventions, we define a dilation-weighted discrete Sobolev space by
\[
\ell^1(w^{m,p},\alpha) = \left\{ c \in \ell^1(\ell^p) :
\text{$\alpha^{|\rho|} \Delta^\rho c \in \ell^1(\ell^p)$ for each
multiindex $\rho$ of order $|\rho| \leq m$} \right\} ,
\]
with norm
\begin{align*}
\| c \|_{\ell^1(w^{m,p},\alpha)} & := \sum_{|\rho| \leq m} \|
\alpha^{|\rho|} \Delta^\rho c \|_{\ell^1(\ell^p)} \\
& = \sum_{j>0} \sum_{|\rho| \leq m} |\alpha_j|^{|\rho|} \left(
\sum_{k \in \Zd} |(\Delta^\rho c)_{j,k}|^p  \right)^{\! \! \! 1/p} .
\end{align*}
(The difference operators here should be understood as acting on the
$k$-index of the sequence $c_{j,k}$.) One can check
$\ell^1(w^{m,p},\alpha)$ is a Banach space.

\vspace{6pt} \emph{The class of synthesizers.} Our synthesis
results will hold when $\psi$ has the special convolution form
\begin{equation} \label{psidef}
\psi = \overset{\text{$m$ factors}}{\overbrace{\beta * \cdots *
\beta}} * \eta ,
\end{equation}
where
\[
\beta=|b\cube|^{-1} \charfn_{b\cube}
\]
is the normalized indicator function of the period box in the
lattice $b\Zd$. (One can show that the convolution form
\eqref{psidef} is equivalent in one dimension to a Strang--Fix
condition on $\psi$, but it is definitely stronger in higher
dimensions, as explained in \cite[\emph{Notes on Theorem~1}]{bl4}.)
The point of this convolution form \eqref{psidef} is that
derivatives of $\psi$ turn into differences, in formula
\eqref{diffformula} later on, and these differences transfer to the
coefficient sequence in formula \eqref{keystep}.

Now we can state the boundedness of the synthesis operator on
Sobolev space. Our statement involves the matrix rescaling operator
\[
(M_b f)(x) = |\det b| f(bx) .
\]
Choosing $b=I$ gives the simplest results in what follows, of
course.
\begin{theorem}[Synthesis into $W^{m,p}$] \label{sobolevrep1}
Assume $1 \leq p \leq \infty$ and $\eta \in L^p$ with $P|\eta| \in
L^p_{loc}$. Let $m \in \N$ and define $\psi$ by \eqref{psidef}.

Then $\psi \in W^{m,p}$, and $S : \ell^1(w^{m,p},\alpha) \to
W^{m,p}$ is bounded. More precisely, if $c \in
\ell^1(w^{m,p},\alpha)$ then the series $Sc = \sum_{j>0} \sum_{k \in
\Zd} c_{j,k} \psi_{j,k}$ converges in $W^{m,p}$ in the sense that
\begin{equation} \label{sobolevconverge}
\begin{minipage}{.9\linewidth}
the sum over $k$ converges pointwise absolutely a.e.\ to a function
in $W^{m,p}$ and the sum over $j$ converges absolutely in
$W^{m,p}$,
\end{minipage}
\end{equation}
and furthermore
\[
\| M_b Sc \|_{W^{m,p}} \leq |\det b|^{-1/p} \| P|\eta|
\|_{L^p(b\cube)} \| c \|_{\ell^1(w^{m,p},\alpha)} .
\]
\end{theorem}
See Section~\ref{sobolevrep1_proof} for the proof. Prior results on
bounded synthesis include \cite{G69,S93}, for compactly supported
$\psi$. Here $\psi$ need not be compactly supported.

Regarding the appearance of the rescaling operator $M_b$ in the
theorem, note the analogous $L^p$ inequality \eqref{analogue1} in
Theorem~\ref{lebesguerep1} can be put in the same form as above,
namely
\[
\| M_b Sc \|_p \leq |\det b|^{-1/p} \| P|\psi| \|_{L^p(b\cube)} \| c
\|_{\ell^1(\ell^p)} ,
\]
by applying Theorem~\ref{lebesguerep1} and then making a change of
variable on the lefthand side.

Now we obtain surjectivity of $S$.
\begin{theorem}[Synthesis onto $W^{m,p}$] \label{sobolevrep2}
Assume $1 \leq p < \infty$ and $\eta \in L^p$ with $P|\eta| \in
L^p_{loc}$ and $\int_\Rd \eta \, dx = 1$. Let $m \in \N$ and define
$\psi$ by \eqref{psidef}.

Then $S : \ell^1(w^{m,p},\alpha) \to W^{m,p}$ is open, and
surjective. Indeed if $f \in W^{m,p}$ and $\e>0$ then a sequence $c
\in \ell^1(w^{m,p},\alpha)$ exists such that $Sc=M_b^{-1} f$ with
convergence as in \eqref{sobolevconverge}, and such that
\[
\| c \|_{\ell^1(w^{m,p},\alpha)} \leq \| f \|_{W^{m,p}} + \e.
\]
\end{theorem}
The theorem is proved in Section~\ref{sobolevrep2_proof}. Note the
conclusion of the $L^p$ result Theorem~\ref{lebesguerep2} can be
rephrased to look like Theorem~\ref{sobolevrep2}, by applying it to
$M_b^{-1} f$ instead of to $f$.

There seem to be no direct predecessors in the literature for the
surjectivity in Theorem~\ref{sobolevrep2}. Indirect predecessors
include the density results from Strang--Fix theory, which we
discuss in \cite[\S3.5]{bl4} and after
Theorem~\ref{scaleaveragesobolev} below. To learn about wavelet and
affine frame expansions in Sobolev space, using all scales $j \in
\Z$, consult \cite{HW96,M92} and the recent work in \cite{CS06}.

\vspace{6pt} Norm equivalence now follows from
Theorems~\ref{sobolevrep1} and \ref{sobolevrep2}:
\begin{corollary}[Synthesis norm for $W^{m,p}$] \label{sobolevnorm1}
Assume $1 \leq p < \infty$ and $\eta \in L^p$ with $P|\eta| \in
L^p_{loc}$ and $\int_\Rd \eta \, dx = 1$. Let $m \in \N$ and define
$\psi$ by \eqref{psidef}.

Then
\[
\| f \|_{W^{m,p}} \approx \inf \left\{ \| c
\|_{\ell^1(w^{m,p},\alpha)} : \text{$M_b^{-1} f = Sc$ \, as in
\eqref{sobolevconverge}} \right\}
\]
for all $f \in W^{m,p}$. Explicitly,
\begin{align*}
|\det b|^{1/p} \| P|\eta| \|_{L^p(b\cube)}^{-1} \| f \|_{W^{m,p}} &
\leq \inf \left\{ \| c \|_{\ell^1(w^{m,p},\alpha)} :
\text{$M_b^{-1} f = Sc$ \, as in \eqref{sobolevconverge}} \right\} \\
& \leq \| f \|_{W^{m,p}} .
\end{align*}
\end{corollary}

The coefficient norm can equal the standard Sobolev norm:
\begin{corollary}[Synthesis norm equality] \label{sobolevnorm2}
Assume $1 \leq p < \infty$ and that $\eta \in L^p$ is nonnegative.
When $p=1$ assume $\int_\Rd \eta \, dx = 1$, and when $1<p<\infty$
assume $P\eta \equiv 1$ (which implies $\int_\Rd \eta \, dx = 1$).
Let $m \in \N$ and define $\psi$ by \eqref{psidef}.

Then for all $f \in W^{m,p}$,
\[
\| f \|_{W^{m,p}} = \inf \left\{ \| c \|_{\ell^1(w^{m,p},\alpha)} :
\text{$M_b^{-1} f = Sc$ \, as in \eqref{sobolevconverge}} \right\} .
\]
\end{corollary}
The analysis operator gives a norm to Sobolev space also:
\begin{corollary}[Analysis norm for
$W^{m,p}$]\label{sobolevnorm3} Assume  $1 \leq p < \infty$ and $m
\in \N$. Take an analyzer $\phi \in L^q$ with $P|\phi| \in L^\infty$
and $\int_{\Rd} \phi \, dx = 1$, and assume $b=I$.

Then for all $f\in W^{m,p}$,
\begin{align*}
\| f \|_{W^{m,p}} \approx \| Tf \|_{\ell^\infty(w^{m,p},\alpha)}  & =
\sup_{j > 0} \sum_{|\rho|\leq m} |\alpha_j|^{|\rho|} \| \Delta^\rho
T_j f \|_{\ell^p}\\
& = \sup_{j > 0} \sum_{|\rho|\leq m} |\alpha_j|^{|\rho|} \|
\{\Delta^\rho \langle f,\phi_{j,k}\rangle\} \|_{\ell^p} .
\end{align*}
\end{corollary}
Thus the analysis operator is linear, bounded and injective from
$W^{m,p}$ onto its range. Section~\ref{sobolevnorm3_proof} has the
proof.

\section{\bf Connection to phi-transforms and Triebel--Lizorkin spaces}
\label{literature}

The decomposition or representation of function spaces by means of
an affine system generated by a single function (or collection of
functions) is a well-established technique in harmonic analysis.
Important examples include wavelet expansions and the
phi-transforms. For wavelet theory we refer to
\cite{Ch03,D92,HW96,M92} and the references therein. The connection
between phi-transform theory and our results is briefly explained
below.

Triebel--Lizorkin spaces include $L^p, H^1$ and $W^{m,p},
1<p<\infty$, but neither $L^1$ nor $W^{m,1}$. Phi transform theory
implies the following results about these spaces; see
\cite{BP04,FJ90,FJW91,GHHLWW} for a complete account.

(i) That the Triebel--Lizorkin space norm of an element $f$ is
equivalent to the infimum of the corresponding sequence space norm
of $c$, where the infimum is taken over all representations $f=Sc$
and the synthesizing $\psi$ satisfies some moment condition (in
particular $\int_\Rd \psi \, dx = 0$), and a decay condition and a
``Tauberian'' condition. Note that the dilations are assumed
isotropic, the translation matrix $b$ must be ``sufficiently
small'', and that the sequence space norm in this theory also
involves an integration with respect to the continuous variable on
$\Rd$.

(ii) That there exists an analyzing function $\phi$ (or a collection
of functions $\phi^{(j,k)}$) satisfying similar conditions to $\psi$
such that the Triebel--Lizorkin space norm of $f$ is equivalent to
the sequence space norm of $Tf = \{\langle f,\phi_{j,k}\rangle\}$,
and $f$ is represented by $f=\sum_{j,k} S_j T_j f$.

The large scales $j \leq 0$ must be included in the synthesis (i)
and the analysis (ii).

Versions of these results have also been proved for affine systems
generated by interesting families of functions (instead of one
function $\psi)$; see \cite{T02} and the references there for
``Gausslet'' and ``Quarkonial'' analysis.

Our main results in Sections
\ref{lebesgueresults}--\ref{sobolevresults} can be viewed as
analogues of the above results when both the analyzer $\phi$ and the
synthesizer $\psi$ have non-zero integral (that is, no vanishing
moments). Further, our sequence space norm is simpler, involving
only discrete sums and no integration with respect to the continuous
variable. And the norm equivalence statements in our
Corollaries~\ref{lebesguenorm3}, \ref{hardynormb} and
\ref{sobolevnorm3} take an especially simple form. These corollaries
are in the spirit of norm equivalences in frame theory, but in our
case the frame decomposition $f=\sum_{j,k} S_j T_j f$ cannot hold,
because if $f \mapsto \sum_{j,k} S_j T_j f$ is bounded on $L^2$,
then either $\phi$ or $\psi$ must have a vanishing moment by
\cite[Theorem B]{GHL05}. Nevertheless, our corollaries give rise to
Banach frames in the sense of K. Gr\"{o}chenig \cite{Gro91}. Details
are given in Appendix~\ref{bframe}.

\section{\bf Open problems}
\label{openproblems}

Throughout this paper we have assumed the synthesizer $\psi$ has
nonvanishing integral. When it has vanishing integral, $\int_\Rd
\psi \, dx = 0$, we do not have a comprehensive understanding of
conditions under which one can synthesize surjectively onto $L^p$ or
Hardy or Sobolev space. Of course there are substantial classes of
synthesizers $\psi$ for which surjectivity and even injectivity is
known, for example the class of wavelets \cite{D92,HW96,M92}
(provided one includes large scales $j \leq 0$ in the affine
systems). But wavelets seem rather restricted objects to us. In
accordance with the goal expressed in the Introduction, we
conjecture that some much more general surjectivity result should
hold when $\psi$ has vanishing integral. Recent $L^2$ work of
Gilbert \emph{et al.} \cite[Theorem~G]{GHL05} is a step in the right
direction, for it assumes only that $\widehat{\psi}(0)=0$ and
$\widehat{\psi}$ has some cancellation properties. Unfortunately the
result suffers from oversampling of both dilations and translations,
and ways to remove that oversampling remain a mystery.

Problems with vanishing integrals can be more challenging than they
appear. For instance, it is an open problem of Y. Meyer
\cite[p.~137]{M92} to determine whether the affine system $\{
\psi(2^j x-k) : j, k \in \Z \}$ spans $L^p(\R)$ for each $1 < p <
\infty$, when $\psi(x)=(1-x^2)e^{-x^2/2}$ is the Mexican hat
function (the second derivative of the Gaussian). This is known to
be true when $p=2$, since the system forms a frame, but it remains
open for all other $p$-values. It is apparently also open to
determine whether the system spans $H^1(\R)$. To express this in
terms of the synthesis operator, first notice $S_\psi$ is bounded
from $\ell^1(\ell^1)$ to $H^1$ because the Mexican hat $\psi$
belongs to $H^1$, and then ask: is $S_\psi : \ell^1(\ell^1) \to
H^1(\R)$ surjective? Surjectivity would give an atomic decomposition
of $H^1$ in terms of the Mexican hat affine system.

There are two partial results in the literature dealing with the
Mexican hat problem. In \cite{GHHLWW} the authors proved that there
exist (sufficiently small) $r>1$ and $s>0$ such that the affine
system $\{\psi(r^j x-sk) : j, k \in \Z \}$ spans $H^p(\R)$, $1\leq
p<\infty$, while in \cite{BP04} the authors proved the same result
for the affine system $\{ \psi(2^j x-bk) : j, k \in \Z \}$, where $b
> 0$ is sufficiently small and $1/2 < p < \infty.$ Since $H^p=L^p$
for $1<p<\infty$, these results show that the spanning property
holds for the Mexican hat function provided we accept some degree of
oversampling.

We conclude by pointing out a gap in understanding of our Hardy
space synthesis results. We do not know whether $h^1$ is the
``natural'' domain for the synthesis operator in
Theorem~\ref{hardyrep1}. Can one prove it is natural, in the sense
that
\[
\sum_{k \in \Zd} s_k \psi_{j,k} \in H^1 \quad \Longrightarrow \quad
s \in h^1
\]
whenever $\psi \in L^1$ satisfies hypothesis \eqref{h1bound}?


\section{\bf $L^p$ proofs}
\label{lebesgue_proofs}

\subsection{\bf Proof of Theorem~\ref{lebesguerep1} ---
$\text{synthesis\ } \ell^1(\ell^p) \to L^p$} \
\label{lebesguerep1_proof}

First assume $1 \leq p < \infty$ and $\psi \in L^p$ with $P|\psi|
\in L^p_{loc}$. We will synthesize at a fixed scale $j>0$, by taking
$s \in \ell^p$ and defining
\[
f = S_j s = \sum_{k \in \Zd} s_k \psi_{j,k} .
\]
The task is to show $f \in L^p$ with $\| f \|_p \leq \| s
\|_{\ell^p} \| P|\psi| \|_{L^p(b\cube)}/|\det b|$. Then
Theorem~\ref{lebesguerep1} follows easily, by summing over the
dilation scales.

We have
\begin{align*}
|\sum_{k \in \Zd} s_k \psi_{j,k}(x)|^p & \leq \left( \sum_{k
\in \Zd} |s_k| |\psi(a_j x - bk)| \right)^{\! \! p} |\det a_j| \\
& \leq \sum_{k \in \Zd} |s_k|^p |\psi(a_j x - bk)| \left( \sum_{k
\in \Zd} |\psi(a_j x - bk)| \right)^{\! p-1} |\det a_j|
\end{align*}
by H\"{o}lder's inequality on the sum. Integrating with respect to
$x$ yields that
\begin{align}
\| \sum_{k \in \Zd} s_k \psi_{j,k} \|_p^p
& \leq \sum_{k \in \Zd} |s_k|^p \int_\Rd |\psi(x - bk)| (P|\psi|(x))^{p-1} \, dx / |\det b|^{p-1} \notag \\
& = \sum_{k \in \Zd} |s_k|^p \| P|\psi| \|_{L^p(b\cube)}^p / |\det
b|^p , \label{kest}
\end{align}
by changing $x \mapsto x+bk$ and then periodizing the integral. We
conclude that the sum over $k$ in \eqref{frep} converges pointwise
absolutely a.e.\ to an $L^p$ function, and that
\begin{align*}
\| f \|_p = \| \sum_{k \in \Zd} s_k \psi_{j,k} \|_p \leq \| s
\|_{\ell^p} \| P|\psi| \|_{L^p(b\cube)}/|\det b| .
\end{align*}

It remains to prove the theorem when $p=\infty$. So assume
$p=\infty$ and $\psi \in L^\infty$ with $P|\psi| \in
L^\infty_{loc}$. Then $P|\psi|$ is bounded, since it is locally
bounded and periodic. If $c \in \ell^1(\ell^\infty)$ then
\begin{align*}
|\sum_{j > 0} \sum_{k \in \Zd} c_{j,k} \psi_{j,k}(x)| & \leq
\sum_{j > 0} (\sup_{k \in \Zd} |c_{j,k}|) \sum_{k \in \Zd}
|\psi(a_j x - bk)| \\
& \leq \| c \|_{\ell^{1,\infty}} \| P|\psi| \|_\infty/|\det b| <
\infty
\end{align*}
for almost every $x$, from which the theorem follows.

\vspace{6pt} \noindent \underline{Is $\ell^1(\ell^p)$ the correct
domain for synthesis?}

We have proved $S$ is bounded from $\ell^1(\ell^p)$ into $L^p$.
Could $S$ be bounded on an even larger domain? The natural candidate
would be $\ell^r(\ell^{p^\prime})$ with $r \geq 1, p^\prime \geq p$,
but we will show by example that $S$ need not be bounded on this
domain unless $r=1$ and $p^\prime=p$.

Work in one dimension with $b=1$ and dyadic dilations $a_j=2^j$, and
choose $\psi$ to be supported in the unit interval $[0,1)$. Then for
any sequence $s \in \ell^{p^\prime}$ we have $S_j s(x) = \sum_{k \in
\Z} s_k 2^{j/p} \psi(2^j x - k)$, which has $L^p$ norm $\| S_j s
\|_p = \| \psi \|_p \| s \|_{\ell^p}$. Hence if $S_j s \in L^p$ then
$s \in \ell^p$, so that for $S$ to be bounded on
$\ell^r(\ell^{p^\prime})$ it is necessary that $p^\prime = p$.

Further, if $t = \{ t_j \}_{j>0}$ is any nonnegative sequence in
$\ell^r$ then the sequence
\[
c_{j,k} =
\begin{cases}
t_j 2^{-j/p} , & k=0,1, \ldots, 2^j - 1 , \\
0 ,  & \text{otherwise,}
\end{cases}
\]
belongs to $\ell^r(\ell^p)$, and if $\psi=\charfn_{[0,1)}$ is the
indicator function of the unit interval then $Sc = \| t \|_{\ell^1}
\charfn_{[0,1)}$. Thus for $Sc$ to belong to $L^p$ it is necessary
that $t \in \ell^1$. Hence $r=1$, as claimed.

So we cannot expect to enlarge the domain $\ell^1(\ell^p)$, in
general, when $\int_\R \psi \, dx \neq 0$. But there is a loophole
relevant to wavelets, because if $\int_\R \psi \, dx = 0$ then the
$\psi_{j,k}$ might exhibit some cancellation between different
$j$-scales that invalidates our ``$t$'' example above and allows us
to take $r>1$. For example, if $p=2$ and $\psi$ is a wavelet
(meaning the functions $\psi_{j,k}$ are orthonormal and complete in
$L^2$) then $S$ is not only bounded but is an isometry from
$\ell^2(\Z \times \Z)$ to $L^2$. Thus the natural domain in the
wavelet case has $r=2$. Recall integrable wavelets satisfy $\int_\R
\psi \, dx = 0$ by \cite[p.~348]{HW96}.

Continuing in the special case $p=2$, the $\psi_{j,k}$ form a frame
(which is more general than an orthonormal basis) if and only if the
synthesis operator is bounded and surjective from $\ell^2(\Z \times
\Z)$ to $L^2$. This is a special case of Christensen's Hilbert space
result \cite[Theorem~5.5.1]{Ch03}. Thus the natural domain for frame
synthesis has $r=2$. This agrees with our earlier remarks, because
integrable frame generators are known to satisfy $\int_\R \psi \, dx
= 0$.

\subsection{\bf $\text{Analysis\ } L^p \to \ell^p$} \
\label{analysislebsec}

The proof of Theorem~\ref{lebesguerep2} relies on the following
estimate for analyzers $\phi$ with bounded periodization. Recall the
analysis operator $T_j$ at scale $j$ was defined in
Section~\ref{definition}.
\begin{proposition}[Analysis into $\ell^p$] \label{analysisleb}
Assume $1 \leq p \leq \infty$ and $\phi \in L^q$ with $P|\phi| \in
L^\infty$. Then for each $j$,
\[
T_j : L^p \to \ell^p \qquad \text{with norm $\| T_j \| \leq |\det
b|^{1/q} \| P|\phi| \|_\infty$.}
\]
\end{proposition}
The proposition is known from Aldroubi, Sun and Tang \cite[formula
(2.2)]{AST01}.

The hypothesis that $|\phi|$ have bounded periodization means $\phi$
is bounded and its integer translates ``do not overlap too often''.
\begin{proof}[Proof of Proposition~\protect\ref{analysisleb}]
Let $f \in L^p$. When $1 \leq p < \infty$ we have
\begin{align*}
& \| T_j f \|_{\ell^p} \\
& = |\det b| \left( \sum_{k \in \Zd} | \langle f , \phi_{j,k}
\rangle |^p \right)^{\! \! \! 1/p} \\
& \leq |\det b| \left( \sum_{k \in \Zd} \int_\Rd |f(y)|^p |\phi(a_j
y - bk)| \, dy \left( \int_\Rd |\phi(a_j y -bk)| |\det a_j| \, dy
\right)^{\! \! p/q} \right)^{\! \! \! 1/p} \\
& \qquad \qquad \qquad \qquad \qquad
\text{by H\"{o}lder's inequality on the inner product} \notag \\
& = |\det b|^{1/q} \left( \int_\Rd |f(y)|^p P|\phi|(a_j y) \, dy
\right)^{\! \! 1/p} \| \phi \|_1^{1/q} \\
& \leq |\det b|^{1/q} \| P|\phi| \|_\infty \| f \|_p ,
\end{align*}
using that $\| \phi \|_1 = \int_{b\cube} \sum_{k \in \Zd}
|\phi(x-bk)| \, dx \leq \| P|\phi| \|_\infty$.

When $p=\infty$, the proof is straightforward.
\end{proof}
\emph{Aside.} For the special case $p=2$,
Proposition~\ref{analysisleb} is known \cite[Theorem~7.2.3]{Ch03}
for all $\phi \in L^2$ with $P(|\widehat{\phi}|^2) \in L^\infty$,
and one can show this condition is weaker than $P|\phi| \in
L^\infty$. That is, analysis is bounded from $L^2$ to $\ell^2$ with
norm $\| T_j \| \leq \| P(|\widehat{\phi}|^2) \|_\infty^{1/2}$,
where the periodization is with respect to the lattice $\Zd b^{-1}$.
This is another way of saying the translates of $\phi$ form a Bessel
sequence, or satisfy an upper frame bound.

\subsection{\bf Scale-averaged approximation in $L^p$} \ \label{scaleaveragesec}

The following approximation result will be used in proving
Theorem~\ref{lebesguerep2} (surjectivity of the synthesis operator).
The result is interesting in its own right too, due to its explicit
nature: we simply analyze with $\phi$, then synthesize with $\psi$,
and then average over all dilation scales to recover $f$.

Say that the dilations \emph{expand exponentially} if
\[
\| a_j a_{j+1}^{-1} \| \leq \delta \qquad \text{for all $j>0$,}
\]
for some $0<\delta<1$. In one dimension, this means $|a_j|$ is a
lacunary sequence.
\begin{theorem}[\protect{\cite[Theorem~1 and Lemma~2]{bl1}}] \label{sample}
Assume $1 \leq p < \infty$ and $\psi \in L^p, P|\psi| \in L^p_{loc},
\phi \in L^q, P|\phi| \in L^\infty, f \in L^p$. Assume $\int_\Rd
\phi \, dx = 1$, and write $\gamma = \int_\Rd \psi \, dx$.

\noindent (a) \emph{[Constant periodization]} If $P \psi = \gamma$
a.e.\ then
\[
S_j T_j f \to \gamma f \qquad \text{in $L^p$ as $j \to \infty$.}
\]

\noindent (b) \emph{[Scale-averaged approximation]} If the dilations
$a_j$ expand exponentially, then
\[
\frac{1}{J} \sum_{j=1}^J S_j T_j f \to \gamma f \qquad \text{in
$L^p$ as $J \to \infty$.}
\]
\end{theorem}
The scale averaged approximation in part (b) can be written in full
as
\[
\frac{1}{J} \sum_{j=1}^J |\det b| \left( \sum_{k \in \Zd} \langle f
, \phi_{j,k} \rangle \psi_{j,k} \right) \to \gamma f  \quad \text{in
$L^p$, as $J \to \infty$.}
\]
\emph{Aside.} Part (a) of the theorem has a long history, summarized
in \cite[\S3]{bl1}.

\subsection{\bf Proof of Theorem~\ref{lebesguerep2} --- synthesis
onto $L^p$} \ \label{lebesguerep2_proof}

We can assume the dilations expand exponentially, as follows. For
each $j$ we have $\| a_j a_{j+r}^{-1} \| \leq \| a_j \| \|
a_{j+r}^{-1} \| \to 0$ as $r \to \infty$ because the dilations are
expanding. Thus $\| a_j a_{j+r}^{-1} \| \leq 1/2$ provided we choose
$r$ sufficiently large. By iterating this argument we arrive at a
subsequence of dilations that expands exponentially (with
$\delta=1/2$). It is enough to use only this subsequence of
dilations, when proving openness and surjectivity of the synthesis
operator.

Consider $f \in L^p$ and $J \in \N$, and define a sequence $c_J = \{
c_{J;j,k} \}_{j>0,k \in \Zd}$ by
\begin{equation} \label{cJdef}
c_{J;j,k} = \frac{1}{J} |\det b|
\begin{cases}
\langle f , \phi_{j,k} \rangle & \text{for $j=1,\ldots,J$,} \\
0 & \text{otherwise,}
\end{cases}
\end{equation}
where $\phi = |b\cube|^{-1} \charfn_{b\cube}$ is a normalized
indicator function. Note $P|\phi| \equiv 1$. Then $c_J \in
\ell^1(\ell^p)$ because applying Proposition~\ref{analysisleb} for
each $j=1,\ldots,J$ gives that
\[
\| c_J \|_{\ell^1(\ell^p)} \leq \frac{1}{J} \sum_{j=1}^J |\det
b|^{1/q} \| f \|_p = |\det b|^{1/q} \| f \|_p .
\]
And clearly
\begin{align*}
Sc_J & = \frac{1}{J} \sum_{j=1}^J S_j T_j f \\
& \to f \qquad \qquad \text{in $L^p$ as $J \to \infty$,}
\end{align*}
by Theorem~\ref{sample}(b). (Here we use that the $a_j$ expand
exponentially.) Thus the open mapping theorem in
Appendix~\ref{banachapp} says $S : \ell^1(\ell^p) \to L^p$ is open
and surjective, and that for each $f \in L^p$ and $\e>0$ there
exists $c \in \ell^1(\ell^p)$ with $Sc = f$ and $\| c
\|_{\ell^1(\ell^p)} \leq |\det b|^{1/q} \| f \|_p + \e$.

\subsection{\bf Proof of Corollary~\ref{lebesguenorm4} ---
Synthesis norm for $L^p(\Omega)$} \ \label{lebesguenorm4_proof}

The ``$\leq$'' direction of the norm equivalence follows immediately
from Corollary~\ref{lebesguenorm1}, since we are restricting the
collection of sequences $c$ that can be used to represent $f$.

For the ``$\geq$'' direction, we will modify the proof of
Theorem~\ref{lebesguerep2}. Like in that theorem we can assume the
dilations expand exponentially, by passing to a subsequence of
$j$-values.

Suppose $f \in L^p(\Omega)$ is supported on a compact subset of
$\Omega$. We claim there exists $j_0 \geq 0$ such that
\begin{equation} \label{suppeq}
\text{$\spt(\psi_{j,k}) \subset \Omega$ whenever $j>j_0$ and
$\langle f , \phi_{j,k} \rangle \neq 0$,}
\end{equation}
where $\phi = |b\cube|^{-1} \charfn_{b\cube}$ as previously. The
existence of $j_0$ should be clear intuitively, since the support of
$f$ lies at some positive distance from the boundary of $\Omega$
while $\psi$ and $\phi$ have compact support and $\| a_j^{-1} \| \to
0$. We leave the detailed proof of \eqref{suppeq} to the reader.

Next, consider $J \in \N$ and define
\[
c_{J;j,k} = \frac{1}{J} |\det b|
\begin{cases}
\langle f , \phi_{j,k} \rangle & \text{for $j=j_0+1,\ldots,j_0+J$,} \\
0 & \text{otherwise,}
\end{cases}
\]
so that $c_J \in \ell^1(\ell^p)$ by Proposition~\ref{analysisleb},
with
\begin{equation} \label{lrep2eq1again}
\| c_J \|_{\ell^1(\ell^p)} \leq \frac{1}{J} \sum_{j=j_0+1}^{j_0+J}
|\det b|^{1/q} \| f \|_p = |\det b|^{1/q} \| f \|_p .
\end{equation}
We know $Sc_J \to f$ in $L^p(\Rd)$ as $J \to \infty$, by
Theorem~\ref{sample}(b) with an index shift on the dilations.

Define
\[
{\mathcal L} = \{ c \in \ell^1(\ell^p) : \text{$c$ is adapted to
$\Omega$ and $\psi$} \} .
\]
Clearly ${\mathcal L}$ is a closed subspace of $\ell^1(\ell^p)$, and
hence is a Banach space. Note $c_J \in {\mathcal L}$ by
\eqref{suppeq}, so that $Sc=0$ on $\Omega^c$.

We have verified the hypotheses of the open mapping theorem in
Appendix~\ref{banachapp} for $S : {\mathcal L} \to L^p(\Omega)$,
with the constant $A=|\det b|^{1/q}$ by \eqref{lrep2eq1again}.
Admittedly we have verified the hypotheses only for the dense class
of $f$ having compact support in $\Omega$, but a dense class
suffices, by the comment at the end of Appendix~\ref{banachapp}. The
open mapping theorem tells us $S : {\mathcal L} \to L^p(\Omega)$ is
surjective, and that for each $f \in L^p(\Omega)$ and $\e>0$ there
exists $c \in {\mathcal L}$ with $Sc = f$ and $\| c
\|_{\ell^1(\ell^p)} \leq |\det b|^{1/q} \| f \|_p + \e$. This proves
the ``$\geq$'' direction of the corollary.

\subsection{\bf Proof of Corollary~\ref{lebesguenorm3} --- analysis norm for $L^p$} \
\label{lebesguenorm3_proof}

By boundedness of the analysis operator in
Proposition~\ref{analysisleb},
\[
\sup_j \| T_j f \|_{\ell^p} \leq |\det b|^{1/q}
\|P|\phi|\|_\infty\|f\|_p.
\]
To prove the other direction of the norm equivalence, choose a
function $\psi$ that satisfies the hypotheses of
Theorem~\ref{sample}(a) with $\gamma=1$. Then by that theorem, $S_j
T_j f \to f$ in $L^p$ as $j\to \infty$. Therefore it follows from
Theorem~\ref{lebesguerep1} (bounded synthesis) that
\[
\| f \|_p \leq \sup_j \| S_j T_j f \|_p \leq |\det b|^{-1}
\|P|\psi|\|_{L^p(b\cube)} \sup_j \| T_j f \|_{\ell^p} ,
\]
which proves the corollary.

\section{\bf Hardy space proofs}
\label{hardy_proofs}

In this section we assume the dilations are isotropic and expanding,
so that $a_j = \alpha_j I$ for some nonzero real numbers $\alpha_j$
with $|\alpha_j| \to \infty$.

\subsection{\bf Proof of Theorem~\ref{hardyrep1} ---
$\text{synthesis\ } \ell^1(h^1) \to H^1$} \ \label{hardyrep1_proof}

We begin by showing that the Riesz transform ``almost'' commutes
with affine synthesis,  which will be the key step in proving
boundedness of the synthesis operator, in Theorem~\ref{hardyrep1}.
\begin{lemma}[Riesz transforms commute with synthesis] \label{rieszsynthesis}
Suppose $\psi \in L^1$. Take $\nu$ to be the smooth, compactly
supported cut-off function used to define $z$ in formula
\eqref{zdefn}. Let $\mu$ be the Schwartz function with
$\widehat{\mu}(\xi) = \nu(\xi b)$, and let $\lambda$ be a Schwartz
function with $\widehat{\lambda}$ supported in $\cube_0 b^{-1}$
and with $\widehat{\lambda}(0)=1$.

If $s \in h^1$ then
\begin{align*}
Z*(\sum_{k \in \Zd} s_k (\psi * \lambda * \mu)(x - bk))
& = \sum_{k \in \Zd} (z*s)_k (\psi * \lambda)(x - bk) \\
& \in L^1 .
\end{align*}
\end{lemma}
We can rephrase the lemma (after rescaling $x \mapsto \alpha_j x$
and using the dilation invariance of the Riesz transform) as
saying
\[
R S_{j,\psi * \lambda * \mu} s = \sign(\alpha_j) S_{j,\psi *
\lambda} R s , \qquad s \in h^1,
\]
where on the lefthand side $R$ denotes the continuous Riesz
transform (convolution with the Riesz kernel $Z$) and on the
righthand side $R$ denotes the discrete Riesz transform
(convolution with the discrete Riesz kernel $z$). Thus we see
Lemma~\ref{rieszsynthesis} is a discrete analogue of the formula
$R(\psi * f)=\psi * (Rf)$, once we remember that affine synthesis
is a discrete analogue of convolution with a synthesizer.
\begin{proof}[Proof of Lemma~\ref{rieszsynthesis}]
Define $f(x)=\sum_{k \in \Zd} s_k (\psi * \lambda * \mu)(x-bk)$.
This sum converges absolutely in $L^1$, because
\[
\| f \|_1 \leq \sum_{k \in \Zd} |s_k| \| \psi * \lambda * \mu(\cdot
- bk) \|_1 = \| s \|_{\ell^1} \| \psi * \lambda * \mu \|_1.
\]
Our task is to show $f \in H^1$, with its Riesz transform being as
stated in the lemma.

Consider the periodic functions
\[
\sigma(\xi) = \sum_{k \in \Zd} s_k e^{-2\pi i \xi k},  \qquad
\zeta(\xi) = \sum_{k \in \Zd} z_k e^{-2\pi i \xi k} ,
\]
which are well defined since $s \in \ell^1$ and $z \in \ell^2$. We
have
\begin{align}
-i \frac{\xi}{|\xi|} \widehat{f}(\xi) & = -i \frac{\xi}{|\xi|}
\sigma(\xi b) \widehat{\psi}(\xi) \widehat{\lambda}(\xi)
\widehat{\mu}(\xi) \notag \\
& = \zeta(\xi b) \sigma(\xi b) \widehat{\psi}(\xi)
\widehat{\lambda}(\xi) \label{trans1}
\end{align}
by definition of the Riesz kernel sequence $z$ in \eqref{zdefn},
using here that $\widehat{\lambda}(\xi)=0$ when $\xi b \notin
\cube_0$.

Of course $\zeta(\xi) \sigma(\xi) = \sum_{k \in \Zd} (z*s)_k
e^{-2\pi i \xi k}$ in $L^2(\cube_0)$, by computing Fourier
coefficients of the two sides in $L^2(\cube_0)$ and using $z \in
\ell^2, s \in \ell^1, z*s \in \ell^2$. Therefore \eqref{trans1} says
\[
-i \frac{\xi}{|\xi|} \widehat{f}(\xi) = \text{Fourier transform of\
} \sum_{k \in \Zd} (z*s)_k (\psi * \lambda)(x-bk) ,
\]
where we note that $z * s \in \ell^1$ by the hypothesis $s \in h^1$.
Thus $(Z * f)(x)=\sum_{k \in \Zd} (z*s)_k (\psi * \lambda)(x-bk) \in
L^1$. This proves the lemma.
\end{proof}

\noindent{\em Proof of Theorem~\ref{hardyrep1}}. Fix $j>0$ and take
$s \in h^1$. Define
\[
f = \sum_{k \in \Zd} s_k \psi_{j,k} .
\]
We will show that $f \in H^1$ with $\| f \|_{H^1} \leq C \| s
\|_{h^1}$, where $C=C(\psi,b)$ is independent of $j$. Then
Theorem~\ref{hardyrep1} follows by summing over $j$.

It is enough to show that the function
\[
g(x) = (M_{a_j}^{-1} f)(x) = \sum_{k \in \Zd} s_k \psi(x-bk)
\]
belongs to $H^1$ with $\| g \|_{H^1} \leq C \| s \|_{h^1}$, because
$\| f \|_{H^1} = \| g \|_{H^1}$ (using that $a_j=\alpha_j I$ is
isotropic). Obviously $\| g \|_1 \leq \| s \|_{\ell^1} \| \psi
\|_1$,  and so our task is to show $\| Rg \|_1 \leq C \| s
\|_{h^1}$.

To understand $Rg$, take $\nu, \mu$ and $\lambda$ as in
Lemma~\ref{rieszsynthesis}, and decompose
\begin{align*}
g & = g_1 + g_2 \qquad \text{where}\\
g_1(x) & = \sum_{k \in \Zd} s_k (\psi * \lambda * \mu)(x-bk) , \\
g_2(x) & = \sum_{k \in \Zd} s_k (\psi - \psi * \lambda * \mu)(x-bk)
.
\end{align*}
Lemma~\ref{rieszsynthesis} implies
\[
\| Rg_1 \|_1 \leq \| z*s \|_{\ell^1} \| \psi * \lambda \|_1 \leq \|
s \|_{h^1} \| \psi * \lambda \|_1 .
\]
(Thus we see Lemma~\ref{rieszsynthesis} is used to push the Riesz
transform onto the \emph{coefficient} sequence $s$, which belongs to
$h^1$, rather than onto the synthesizer $\psi * \lambda * \mu$,
which does not belong to $H^1$.)

The sum defining $g_2$ converges absolutely in $H^1$, because $s \in
h^1 \subset \ell^1$ by assumption and $\psi - \psi * \lambda * \mu
\in H^1$ by Lemma~\ref{smoothing} below (which is where we use the
hypothesis \eqref{h1bound}). Hence
\[
\| Rg_2 \|_1 \leq \| s \|_{\ell^1} \| R(\psi - \psi * \lambda * \mu)
\|_1 .
\]
These bounds prove Theorem~\ref{hardyrep1}.

\vspace{6pt} We must still prove Lemma~\ref{smoothing}, needed to
treat $g_2$ in the ``smoothing'' step of the proof above. We use the
translation operator $\tau_y \psi = \psi(\cdot - y)$.
\begin{lemma} \label{smoothing}
Suppose $\psi \in L^1$ and $\| R( \psi - \tau_y \psi ) \|_1 \leq C
< \infty$ for all $|y| \leq 1$, and $\lambda \in L^1$ with
$|y|\lambda(y) \in L^1$ and $\int_\Rd \lambda(y) \, dy =1$. Then
\[
\| R(\psi -  \psi * \lambda) \|_1 \leq C \int_\Rd (|y|+1)
|\lambda(y)| \, dy .
\]
\end{lemma}
The assumption that $\| R( \psi - \tau_y \psi ) \|_1 \leq C$ for
$|y| \leq 1$ is a restatement of hypothesis \eqref{h1bound}.
\begin{proof}[Proof of Lemma~\ref{smoothing}]
The first step is to show the Hardy norm of a difference grows at
most linearly with the difference step, that is
\begin{equation} \label{lingrowth}
\| R( \psi - \tau_y \psi ) \|_{L^1} \leq C(|y|+1) \qquad \text{for
all $y \in \Rd$.}
\end{equation}
For this, let $m$ be the integer satisfying $|y| < m \leq |y|+1$.
After writing $y$ as a sum of $m$ vectors each having norm less
than $1$, we can prove \eqref{lingrowth} with $Cm$ on the
righthand side by telescoping the differences and using the
triangle inequality, noting the Riesz transform is translation
invariant.

Now observe the function
\[
\psi(x) - (\psi * \lambda)(x) = \int_\Rd (\psi(x)-\tau_y \psi(x))
\lambda(y) \, dy
\]
belongs to $H^1$ and has Riesz transform
\[
R(\psi - \psi * \lambda)(x) = \int_\Rd R(\psi-\tau_y \psi)(x)
\lambda(y) \, dy ,
\]
by \cite[Lemma~10]{bl3}. That is, one can take the Riesz transform
through the integral. Thus
\begin{align*}
\| R( \psi - \psi * \lambda) \|_1 & \leq \int_\Rd \| R(\psi - \tau_y \psi) \|_1 |\lambda(y)| \, dy \\
& \leq \int_\Rd C(|y|+1) |\lambda(y)| \, dy
\end{align*}
by \eqref{lingrowth}, as desired.
\end{proof}

\subsection{\bf $\text{Analysis\ } H^1 \to h^1$} \
\label{analysishardysec}

For the proof of Theorem~\ref{hardyrep2} we want boundedness of the
analysis operator $T_j$ from $H^1$ to $h^1$, which we state in
Proposition~\ref{analysishardy} below. First we show:
\begin{lemma}[Riesz transforms commute with analysis] \label{rieszanalysis}
Assume the analyzer $\phi$ is a Schwartz function with
$\widehat{\phi}$ supported in $\cube_0 b^{-1}$. Fix $j>0$. Choose
$\nu$ as in the definition of $h^1$ in Section~\ref{hardyresults},
and let $\mu_j$ be the Schwartz function with $\widehat{\mu_j}(\xi)
= \nu(\xi a_j^{-1} b)$.

If $f \in H^1$ then
\[
R T_j f = \sign(\alpha_j) T_j R(\mu_j * f)  ,
\]
where on the lefthand side $R$ denotes the discrete Riesz transform
(convolution with the discrete Riesz kernel $z$) and on the
righthand side $R$ denotes the continuous Riesz transform
(convolution with the Riesz kernel $Z$).
\end{lemma}
Remember $\alpha_j$ denotes the isotropic dilation factor in
$a_j=\alpha_j I$.
\begin{proof}[Proof of Lemma~\ref{rieszanalysis}]
Observe $\phi_{j,k}(x) = \phi(a_j x - bk)$ in what follows, because
we implicitly assume $p=1,q=\infty$, wherever we deal with the Hardy
space.

The $k$-th term of the sequence $R T_j f$ is
\begin{align*}
(z*T_j f)_k
& = |\det b| \sum_{\ell \in \Zd} z_\ell \langle f , \phi_{j,k-\ell} \rangle \\
& = |\det b| \sum_{\ell \in \Zd} \int_{\cube_0 b^{-1}} z_\ell
e^{-2\pi i \xi b \ell} \widehat{f}(\xi a_j)
\overline{\widehat{\phi}(\xi)e^{-2\pi i \xi bk}} \, d\xi ,
\end{align*}
by Plancherel and the compact support of $\widehat{\phi}$.
Substituting in the definition of $\zeta$ from \eqref{zdefn} and
then using the definition of $\mu_j$, we find
\begin{align*}
& (z*T_j f)_k \\
& = |\det b| \int_\Rd \frac{-i\xi}{|\xi|} \nu(\xi b)
\widehat{f}(\xi a_j) \overline{\widehat{\phi}(\xi) e^{-2\pi i \xi bk}} \, d\xi \\
& = |\det b| \int_\Rd \sign(\alpha_j) (\mu_j * Rf)\widehat{\ \,
}(\xi a_j) \overline{\widehat{\phi}(\xi) e^{-2\pi i \xi bk}} \, d\xi
\qquad \text{since $a_j = \alpha_j I$}
\\
& = |\det b| \sign(\alpha_j) \langle R(\mu_j * f) , \phi_{j,k}
\rangle
\end{align*}
by Plancherel, and this is the $k$-th term of $\sign(\alpha_j) T_j
R(\mu_j * f)$, as desired.
\end{proof}
\begin{proposition}[Analysis into $h^1$] \label{analysishardy}
Take $\phi$ and $\nu$ as in Lemma~\ref{rieszanalysis}. Then for each
$j$,
\[
T_j : H^1 \to h^1 \qquad \text{with norm $\| T_j \| \leq \| P|\phi|
\|_\infty \| \widechecknu \|_1$.}
\]
\end{proposition}
\begin{proof}[Proof of Proposition~\protect\ref{analysishardy}]
If $f \in H^1$ then
\[
\| T_j f \|_{\ell^1} \leq \| P|\phi| \|_\infty \| f \|_1
\]
by Proposition~\ref{analysisleb} with $p=1$,
while
\begin{align*}
\| z * T_j f \|_{\ell^1} = \| R T_j f \|_{\ell^1}
& \leq \| P|\phi| \|_\infty \| R(\mu_j * f) \|_1 \\
& \leq \| P|\phi| \|_\infty \| \mu_j \|_1 \| Rf \|_1
\end{align*}
by combining Lemma~\ref{rieszanalysis} and
Proposition~\ref{analysisleb}. Add these two estimates and observe
$\| \mu_j \|_1 = \| \widechecknu \|_1 \geq |\nu(0)|=1$, by
definition of $\widehat{\mu_j}$ in Lemma~\ref{rieszanalysis}.
\end{proof}

\vspace{6pt} \noindent \emph{Aside.} Our compact support
assumption on the Fourier transform of the analyzer $\phi$, in
Lemma~\ref{rieszanalysis} and Proposition~\ref{analysishardy},
seems rather strong. Perhaps it can be weakened. But notice it
ensures that $\widehat{\phi}(\ell b^{-1})=0$ for all row vectors
$\ell \in \Zd \setminus \{ 0 \}$, which implies $P\phi \equiv
\text{const.}$ This constant periodization condition is necessary,
as follows. If $T_j$ maps $H^1$ into $h^1$, then for all $f \in
H^1$ we have
\[
\int_\Rd f(x) \overline{P\phi(a_jx)} \, dx = \sum_{k \in \Zd} (T_j
f)_k = 0
\]
by the zero-mean property of $T_j f \in h^1$ (see
Appendix~\ref{hardyapp}). Taking $f \in H^1$ to approach a
difference of delta functions implies $P\phi$ is constant.

\subsection{\bf Scale-averaged approximation in $H^1$} \
\label{scaleaveragehardysec}

In this section we prove scale-averaged convergence in $H^1$, which
we need for the proof of Theorem~\ref{hardyrep2}.
\begin{theorem} \label{scaleaveragehardy}
Assume $\psi \in L^1$ with $\int_\Rd \psi \, dx = 1$ and
\begin{equation} \label{hlimitagain}
\| \psi - \psi(\cdot - y) \|_{H^1} \to 0 \qquad \text{as $y \to 0$,}
\end{equation}
which is hypothesis \eqref{h1limit}. Assume $\phi$ is a Schwartz
function with $\widehat{\phi}$ supported in $\cube_0 b^{-1}$ and
$\int_\Rd \phi \, dx = 1$. Let $f \in H^1$.

\noindent (a) \emph{[Constant periodization]} If $P \psi = 1$ a.e.\
then
\[
S_j T_j f \to f \qquad \text{in $H^1$ as $j \to \infty$.}
\]

\noindent (b) \emph{[Scale-averaged approximation]} If the dilations
$a_j$ expand exponentially, then
\[
\frac{1}{J} \sum_{j=1}^J S_j T_j f \to f \qquad \text{in $H^1$ as $J
\to \infty$.}
\]
\end{theorem}
The theorem was proved in our earlier paper \cite[Theorem~1]{bl3}
for analyzers $\phi$ with compact support and $P\phi \equiv
\text{const.}$, by comparing $S_j T_j f$ with an approximate
identity formula. The new proof below is more conceptually
satisfying, as it is based on commuting the Riesz transform
through the analysis and synthesis operators.
\begin{proof}[Proof of Theorem~\ref{scaleaveragehardy}]
Convergence in $L^1$, for parts (a) and (b), follows immediately
from our $L^p$ result Theorem~\ref{sample}, with $p=1$.

To prove convergence in $H^1$, we want to show
\begin{align}
R S_j T_j f \to Rf & \quad \text{in $L^1$, for part (a), \ and} \label{desired1} \\
\frac{1}{J} \sum_{j=1}^J R S_j T_j f \to Rf & \quad \text{in $L^1$,
for part (b).} \label{desired2}
\end{align}
To begin, suppose $\widehat{\mu}(\xi)=\nu(\xi b)$ and $\lambda$ are
as in Lemma~\ref{rieszsynthesis}, and decompose
\begin{align*}
S_j T_j f = S_{j,\psi} T_j f & = S_{j,\psi - \psi * \lambda * \mu}
T_j f + S_{j,\psi * \lambda * \mu} T_j f
\\
& = A_j + B_j , \quad \text{say,}
\end{align*}
where we write $S_{j,\psi}$ and so on to emphasize the synthesizer
being used, in each part of the formula. We have
\begin{align*}
RB_j & = \sign(\alpha_j) S_{j,\psi * \lambda} R T_j f \qquad \text{by Lemma~\ref{rieszsynthesis}} \\
& = S_{j,\psi * \lambda} T_j R(\mu_j * f) \qquad \text{by Lemma~\ref{rieszanalysis}} \\
& = S_{j,\psi * \lambda} T_j ( \mu_j * (Rf) - Rf ) + S_{j,\psi * \lambda} T_j R f \\
& = C_j + D_j , \quad \text{say.}
\end{align*}

We estimate $C_j$ by
\begin{align*}
\| C_j \|_1 & \leq (\text{const.}) \| \mu_j * (Rf) -
Rf \|_1 \qquad \text{by Theorem~\ref{lebesguerep1} and Proposition~\ref{analysisleb}} \\
& \to 0 \qquad \text{as $j \to \infty$,}
\end{align*}
since $\mu_j(x) = |\det b^{-1} a_j| \widechecknu(b^{-1} a_j x)$ is
an approximation to the identity (recalling $a_j=\alpha_j I$ with
$|\alpha_j| \to \infty$).

Theorem~\ref{sample} implies $D_j \to Rf$ in $L^1$ in part (a), and
implies in part (b) that $\frac{1}{J} \sum_{j=1}^J D_j \to Rf$ in
$L^1$.

Thus to prove \eqref{desired1}--\eqref{desired2}, it suffices to
show $RA_j \to 0$ in part (a), and that $\frac{1}{J} \sum_{j=1}^J
RA_j \to 0$ in part (b). To accomplish this, first compute
\begin{align}
RA_j & = R S_{j,\psi - \psi * \lambda * \mu} T_j f \notag \\
& = \sign(\alpha_j) S_{j,R(\psi - \psi * \lambda * \mu)} T_j f ,
\label{rb}
\end{align}
where it is permissible here to pass the Riesz transform through
the synthesis operator because the series for $S_{j,\psi - \psi *
\lambda * \mu} T_j f$ converges absolutely in $H^1$: the
coefficient sequence $T_j f$ belongs to $\ell^1$ by
Proposition~\ref{analysisleb}, and the synthesizer $\psi - \psi *
\lambda * \mu$ belongs to $H^1$ by Lemma~\ref{smoothing} (noting
\eqref{hlimitagain} implies \eqref{h1bound} by \cite[\S3.4]{bl3}).

Next notice $\int_\Rd R(\psi - \psi * \lambda * \mu) \,  dx = 0$,
since all Hardy space functions and their Riesz transforms integrate
to zero; cf.\ \eqref{intzero}. Thus in part (b) of the theorem we
deduce that $\frac{1}{J} \sum_{j=1}^J RA_j \to 0$, from \eqref{rb}
and Theorem~\ref{sample}(b) (and also splitting the sum
$\sum_{j=1}^J$ into two pieces, where $\alpha_j>0$ and $\alpha_j<0$
respectively).

In part (a) of the theorem we deduce that $RA_j \to 0$, by using
\eqref{rb} and Theorem~\ref{sample}(a), and the following
observation. If $P\psi = 1$ a.e., then by computing the Fourier
coefficients of $P\psi$ we find $\widehat{\psi}(\ell b^{-1})=0$ for
all $\ell \in \Zd \setminus \{ 0 \}$, and thus
\begin{equation} \label{intvanish}
(R(\psi - \psi * \lambda * \mu))\widehat{\ \,}(\ell b^{-1}) = 0 .
\end{equation}
Of course \eqref{intvanish} holds for $\ell=0$ too, as observed in
the preceding paragraph. Hence $PR(\psi - \psi * \lambda * \mu) = 0$
a.e., by computing the Fourier coefficients of this periodic
function.

This finishes the proof.
\end{proof}

\subsection{\bf Proof of Theorem~\ref{hardyrep2} --- synthesis onto
$H^1$} \ \label{hardyrep2_proof}

We can assume the dilations expand exponentially, like we did in the
proof of Theorem~\ref{lebesguerep2}.

Let $f \in H^1$ and $J \in \N$, and define the sequence $c_J$ by
\eqref{cJdef} where now $\phi$ is assumed to be  a Schwartz function
with $\widehat{\phi}$ supported in $\cube_0 b^{-1}$ and $\int_\Rd
\phi \, dx = 1$. Then $c_J \in \ell^1(h^1)$ by applying
Proposition~\ref{analysishardy} for each $j=1,\ldots,J$, giving
\begin{equation} \label{hrep2eq1}
\| c_J \|_{\ell^1(h^1)} \leq \frac{1}{J} \sum_{j=1}^J \| P|\phi|
\|_\infty \| \widechecknu \|_1 \| f \|_{H^1} = \| P|\phi|
\|_\infty \| \widechecknu \|_1 \| f \|_{H^1} .
\end{equation}

Observe $Sc_J = \frac{1}{J} \sum_{j=1}^J S_j T_j f \to f$ in $H^1$
as $J \to \infty$ by Theorem~\ref{scaleaveragehardy}. (This is where
the hypotheses on $\psi$ are needed, and that the dilations expand
exponentially.) Combining this with \eqref{hrep2eq1}, we see the
open mapping theorem in Appendix~\ref{banachapp} implies $S$ is
open, and that for each $\e>0$ there exists $c \in \ell^1(h^1)$ with
$Sc = f$ and $\| c \|_{\ell^1(h^1)} \leq \| P|\phi| \|_\infty \|
\widechecknu \|_1 \| f \|_{H^1} + \e$. Choosing $C= \| P|\phi|
\|_\infty \| \widechecknu \|_1$ proves the theorem.

\subsection{\bf Proof of Corollary~\ref{hardynormb} --- analysis norm for $H^1$}\
\label{hardynormb_proof}

Take $\nu$ as in the definition of $h^1$ in
Section~\ref{hardyresults}. Then by boundedness of the analysis
operator in Proposition~\ref{analysishardy},
\[
\sup_j \| T_j f \|_{h^1} \leq \| P|\phi| \|_\infty \| \widechecknu
\|_1 \| f \|_{H^1} .
\]

To prove the other direction of the equivalence, choose $\psi$ to be
a Schwartz function with $P\psi \equiv 1$. Then $S_j T_j f \to f$ in
$H^1$ as $j \to \infty$ by Theorem~\ref{scaleaveragehardy}(a),
noting that hypothesis \eqref{hlimitagain} is known to hold for the
Schwartz function $\psi$ (cf.\ \cite[\S3.3]{bl3}). Therefore it
follows from Theorem~\ref{hardyrep1} (bounded synthesis) that
\[
\|f\|_{H^1} \leq \sup_j \| S_j T_j f \|_{H^1} \leq C \sup_j \| T_j f
\|_{h^1},
\]
where $C$ is independent of $f$.


\section{\bf Sobolev space proofs}
\label{sobolev_proofs}

Throughout this section we assume the dilations are isotropic and
expanding, meaning $a_j = \alpha_j I$ for some nonzero real numbers
$\alpha_j$ with $|\alpha_j| \to \infty$.

We introduce the notation $g^{*r}$ for the convolution of a
function $g$ with itself $r$ times (for example, $g^{*2}=g*g$),
and we write
\[ \Delta_y g = g - g(\cdot - y) , \qquad y \in \Rd,
\]
for the backwards difference of $g$ by $y$. Also we define
\[
\beta_1(x) = \delta(x_1) \charfn_{[0,1)}(x_2) \cdots
\charfn_{[0,1)}(x_d)
\]
and similarly define functions $\beta_2,\ldots,\beta_d$, so that
the partial derivatives of the unit cube indicator function are
\[
D_t \charfn_\cube = \Delta_{e_t} \beta_t, \qquad t=1,\ldots,d.
\]

\subsection{\bf Proof of Theorem~\ref{sobolevrep1} ---
$\text{synthesis\ } \ell^1(w^{m,p},\alpha) \to W^{m,p}$} \
\label{sobolevrep1_proof}

First assume $b=I$, so that $\beta=\charfn_\cube$. Fix $c \in
\ell^1(w^{m,p},\alpha)$.

The initial task is to show that $Sc \in W^{m,p}$. By
differentiating the definition \eqref{psidef} of $\psi$ we convert
derivatives to differences:
\begin{align}
D^\rho \psi & = (D_1 \charfn_\cube)^{*\rho_1} * \cdots
* (D_d \charfn_\cube)^{*\rho_d} * \beta^{* \, m-|\rho|} * \eta \notag \\
& = \Delta_{e_1}^{\rho_1} \cdots \Delta_{e_d}^{\rho_d} (\beta_1^{*
\rho_1} * \cdots * \beta_d^{* \rho_d} * \beta^{* \, m-|\rho|} *
\eta) \in L^p . \label{diffformula}
\end{align}

Obviously \eqref{diffformula} implies that $\psi \in W^{m,p}$.
Further, since $P|\eta| \in L^p_{loc}$ by hypothesis, we have
\begin{align}
P|\beta_1^{* \rho_1} * \cdots * \beta_d^{* \rho_d} * \beta^{* \,
m-|\rho|} * \eta| & \leq \beta_1^{* \rho_1} * \cdots * \beta_d^{*
\rho_d} * \beta^{* \, m-|\rho|} * P|\eta| \label{periodeta} \\
& \in L^p_{loc} , \notag
\end{align}
and therefore $P|D^\rho \psi| \in L^p_{loc}$ by
\eqref{diffformula}. This allows us to use $D^\rho \psi$ as a
synthesizer when applying Theorem~\ref{lebesguerep1}, below.

It is straightforward to show that $Sc$ (which belongs to $L^p$ by
Theorem~\ref{lebesguerep1}) has weak derivatives given by
``differentiating through the sum'', namely
\begin{equation} \label{derivsynth}
D^\rho (S_\psi c) = S_{D^\rho \psi}(\alpha^{|\rho|} c) .
\end{equation}
Note that the righthand side belongs to $L^p$ by
Theorem~\ref{lebesguerep1}, since $c \in \ell^1(w^{m,p},\alpha)$
ensures $\alpha^{|\rho|} c \in \ell^1(\ell^p)$. Thus the function
$Sc$ belongs to $W^{m,p}$, completing the first task in the proof.

The next task is to prove that
\[
\text{\emph{derivatives commute with synthesis,}}
\]
in the sense that
\begin{equation} \label{derivcommute}
D^\rho S_\psi c =  S_{\eta_\rho} (\alpha^{|\rho|} \Delta^\rho c) ,
\end{equation}
where we have introduced the function
\[
\eta_\rho = \beta_1^{* \rho_1} * \cdots * \beta_d^{* \rho_d} *
\beta^{* \, m-|\rho|} * \eta
\]
(so that for example, $\eta_0 = \psi$). Indeed
\begin{align}
(D^\rho S_\psi c)(x)
& = \sum_{j>0} \sum_{k \in \Zd} \alpha_j^{|\rho|} c_{j,k} (D^\rho \psi)_{j,k}(x)
\quad \text{by \eqref{derivsynth}} \notag \\
& = \sum_{j>0} \alpha_j^{|\rho|} \sum_{k \in \Zd} c_{j,k} |\det
a_j|^{1/p} (\Delta_{e_1}^{\rho_1} \cdots \Delta_{e_d}^{\rho_d}
\eta_\rho)(a_j x - k) \quad \text{by \eqref{diffformula}} \notag \\
& = \sum_{j>0} \alpha_j^{|\rho|} \sum_{k \in \Zd} (\Delta^\rho
c)_{j,k} |\det a_j|^{1/p} \eta_\rho(a_j x - k) \label{keystep} \\
& \qquad \qquad \qquad \text{by summation by parts, the key step in the proof,} \notag \\
& = S_{\eta_\rho} (\alpha^{|\rho|} \Delta^\rho c)(x) , \notag
\end{align}
which proves \eqref{derivcommute}.

We will deduce an estimate on the synthesis operator of the form
\begin{equation} \label{middleest}
\| Sc \|_{W^{m,p}} \leq \| P|\eta| \|_{L^p(\cube)} \| c
\|_{\ell^1(w^{m,p},\alpha)} ,
\end{equation}
which completes the proof of the theorem when $b=I$. Start by
observing
\begin{align*}
\| Sc \|_{W^{m,p}} & = \sum_{|\rho| \leq
m} \| D^\rho Sc \|_p \\
& = \sum_{|\rho| \leq m} \| S_{\eta_\rho} (\alpha^{|\rho|}
\Delta^\rho c) \|_p \qquad \text{by \eqref{derivcommute}} \\
& \leq \sum_{|\rho| \leq m} \| P|\eta_\rho| \|_{L^p(\cube)} \|
\alpha^{|\rho|} \Delta^\rho c \|_{\ell^1(\ell^p)}
\end{align*}
by Theorem~\ref{lebesguerep1} (boundedness of synthesis into
$L^p$). To complete the proof of \eqref{middleest}, notice
\[
\| P|\eta_\rho| \|_{L^p(\cube)} \leq \| P|\eta| \|_{L^p(\cube)}
\]
by \eqref{periodeta}.

To prove the theorem when $b \neq I$, first rescale the definition
\eqref{psidef} of $\psi$ to obtain that
\[
M_b \psi = \overset{\text{$m$ factors}}{\overbrace{\charfn_\cube *
\cdots * \charfn_\cube}} * M_b \eta .
\]
That is, \eqref{psidef} holds with $M_b \psi, M_b \eta$ and $I$
instead of $\psi, \eta$ and $b$, so that $M_b \psi$ and $M_b \eta$
satisfy the hypotheses of the theorem for ``$b=I$''.

By the ``$b=I$'' case of the theorem already proved, then, we have
$M_b \psi \in W^{m,p}$ (so that $\psi \in W^{m,p}$), and for each
sequence $c \in \ell^1(w^{m,p},\alpha)$ we have $S_{M_b \psi,I} c
\in W^{m,p}$ with norm estimate
\[
\| S_{M_b \psi,I} c \|_{W^{m,p}} \leq \| P_I|M_b \eta|
\|_{L^p(\cube)} \| c \|_{\ell^1(w^{m,p},\alpha)} .
\]
Further,
\begin{align}
(S_{M_b \psi,I} c)(x) & = \sum_{j>0} \sum_{k \in \Zd} c_{j,k} |\det a_j|^{1/p} (M_b \psi)(a_j x - k) \notag \\
& = \sum_{j>0} \sum_{k \in \Zd} c_{j,k} |\det a_j|^{1/p} |\det b| \psi(b(a_j x - k)) \notag \\
& = |\det b| \sum_{j>0} \sum_{k \in \Zd} c_{j,k} |\det a_j|^{1/p}
\psi(a_j bx - bk) \notag \\
& \qquad \qquad \qquad \text{(noting $b$ commutes with
$a_j=\alpha_j I$)} \notag \\
& = (M_b S_{\psi,b} c)(x) . \label{rescaling}
\end{align}
Hence
\begin{align*}
\| M_b S_{\psi,b} c \|_{W^{m,p}}
& = \| S_{M_b \psi,I} c \|_{W^{m,p}} \\
& \leq \| P_I|M_b \eta| \|_{L^p(\cube)} \| c \|_{\ell^1(w^{m,p},\alpha)} \\
& = |\det b|^{-1/p} \| P|\eta| \|_{L^p(b\cube)} \| c
\|_{\ell^1(w^{m,p},\alpha)} ,
\end{align*}
which finishes the proof.

\subsection{\bf $\text{Analysis\ } W^{m,p} \to w^{m,p}$} \
\label{analysissobsec}

The proof of Theorem~\ref{sobolevrep2} relies on boundedness of
analyzers acting on Sobolev space, as developed in the next
proposition. For simplicity we assume $b=I$, so that the analysis
operator at scale $j$ (defined in Section~\ref{definition}) is just
$(T_j f)_k = \langle f , \phi_{j,k} \rangle$.
\begin{proposition}[Analysis into $w^{m,p}$] \label{analysissob}
Assume $1 \leq p \leq \infty$ and $m \in \N$. Take $\phi \in L^q$
with $P|\phi| \in L^\infty$. Fix $j>0$, and assume $b=I$.

Then $T_j : W^{m,p} \to w^{m,p}$, with norm controlled by the
estimate
\[
\sum_{|\rho| \leq m} |\alpha_j|^{|\rho|} \| \Delta^\rho T_j f
\|_{\ell^p} \leq \| P|\phi| \|_\infty \| f \|_{W^{m,p}} , \qquad f
\in W^{m,p} .
\]
\end{proposition}
\begin{proof}[Proof of Proposition~\protect\ref{analysissob}]
Let $f \in W^{m,p}$ and observe
\[
(T_j f)_k = \langle f , \phi_{j,k} \rangle = \int_\Rd
f(y+\alpha_j^{-1} k) |\det a_j|^{1/q} \overline{\phi(a_j y)} \, dy .
\]
Hence for each multiindex $\rho$ of order $\leq m$ we have
\begin{align*}
\alpha_j^{|\rho|} \Delta^\rho (T_j f)_k & = \int_\Rd
\alpha_j^{|\rho|} \Delta^\rho f(y+\alpha_j^{-1} k) |\det a_j|^{1/q} \overline{\phi(a_j y)} \, dy \\
& = \int_\Rd (\Delta^\rho_{\alpha_j} f)(y+\alpha_j^{-1} k) |\det
a_j|^{1/q} \overline{\phi(a_j y)} \, dy
\end{align*}
where the function $\Delta^\rho_{\alpha_j} f$ in this last line
denotes the $\rho$-th backwards difference quotient of $f$ with step
size $\alpha_j^{-1}$. Changing variable with $y \mapsto y -
\alpha_j^{-1}k$ gives
\begin{equation} \label{diffanal}
\alpha_j^{|\rho|} \Delta^\rho (T_j f)_k = \langle
\Delta^\rho_{\alpha_j} f , \phi_{j,k} \rangle = (T_j
\Delta^\rho_{\alpha_j} f)_k ,
\end{equation}
which says that
\[
\text{\emph{differences commute with the analysis operator.}}
\]
Thus for each fixed $j$, taking the $\ell^p$-norm with respect to
$k$ in \eqref{diffanal} implies
\begin{align*}
\| \alpha_j^{|\rho|} \Delta^\rho T_j f \|_{\ell^p} & \leq \| P|\phi|
\|_\infty \| \Delta^\rho_{\alpha_j} f \|_p \qquad \text{by
Proposition~\ref{analysisleb}} \\
& \leq \| P|\phi| \|_\infty \| D^\rho f \|_p ,
\end{align*}
by using the fundamental theorem of calculus. Summing over $|\rho|
\leq m$ now proves the proposition.
\end{proof}

\subsection{\bf Scale-averaged approximation in $W^{m,p}$} \
\label{scaleaveragesobsec}

Here we prove scale-averaged convergence in $W^{m,p}$, which we
need for the proof of Theorem~\ref{sobolevrep2}. Just like for
$L^p$ and the Hardy space, the idea is to analyze with $\phi$,
then synthesize with $\psi$, and then average over all dilation
scales.
\begin{theorem} \label{scaleaveragesobolev}
Assume $1 \leq p < \infty$ and $\eta \in L^p$ with $P|\eta| \in
L^p_{loc}$ and $\int_\Rd \eta \, dx = 1$. Let $m \in \N$ and define
$\psi$ by \eqref{psidef}. Take $\phi \in L^q$ with $P|\phi| \in
L^\infty$ and $\int_\Rd \phi \, dx = 1$. Assume $b=I$. Let $f \in
W^{m,p}$.

\noindent (a) \emph{[Constant periodization]} If $P \eta = 1$ a.e.\
then
\[
S_j T_j f \to f \qquad \text{in $W^{m,p}$ as $j \to \infty$.}
\]

\noindent (b) \emph{[Scale-averaged approximation]} If the dilations
$a_j$ expand exponentially, then
\[
\frac{1}{J} \sum_{j=1}^J S_j T_j f \to f \qquad \text{in $W^{m,p}$
as $J \to \infty$.}
\]
\end{theorem}
The theorem was proved already in our paper \cite[Theorem~1]{bl4},
by comparing $S_j T_j f$ with an approximate identity formula. The
proof we give below is considerably easier, and is based on passing
derivatives and differences through the analysis and synthesis
operators. On the other hand, our assumptions on the synthesizer
$\psi$ are noticeably stronger than the Strang--Fix type assumptions
in \cite{bl4}, because here we assume $\psi$ has the special
convolution form \eqref{psidef}. This is explained in more detail in
\cite[\emph{Notes on Theorem~1}]{bl4}.

Theorem~\ref{scaleaveragesobolev}(a) is due to Di Guglielmo
\cite[Th\'{e}or\`{e}me~$2^\prime$]{G69}, when $\eta$ has compact
support. Further in this direction, Strang--Fix theory establishes
approximation rates of the firm $O(|\alpha_j|^{k-m})$ in $W^{k,p}$,
for $k<m$, which improves on the rate $o(1)$ in
Theorem~\ref{scaleaveragesobolev}(a). See Jia
\cite[Theorem~3.1]{J04}, and our discussion of the literature in
\cite[\S3.5]{bl4} (where some results special to $p=2$ are cited
also).

Theorem~\ref{scaleaveragesobolev}(b) is needed below when proving
Theorem~\ref{sobolevrep2}.
\begin{proof}[Proof of Theorem~\ref{scaleaveragesobolev}]
Let $\rho$ be a multiindex of order $\leq m$. Then
\begin{align*}
D^\rho S_j T_j f & = S_{j,\eta_\rho} (\alpha^{|\rho|}
\Delta^\rho T_j f) && \text{by \eqref{derivcommute}} \\
& = S_{j,\eta_\rho} T_j \Delta^\rho_{\alpha_j} f && \text{by
\eqref{diffanal}} \\
& = S_{j,\eta_\rho} T_j (\Delta^\rho_{\alpha_j} f - D^\rho f)
+ S_{j,\eta_\rho} T_j D^\rho f \\
& = A_j + B_j , \qquad \text{say.}
\end{align*}
For $A_j$ we have
\begin{align*}
\| A_j \|_p & \leq \| P|\eta_\rho| \|_{L^p(\cube)} \| T_j
(\Delta^\rho_{\alpha_j} f - D^\rho f) \|_{\ell^p}
&& \text{by Theorem~\ref{lebesguerep1}} \\
& \leq \| P|\eta_\rho| \|_{L^p(\cube)} \| P|\phi| \|_\infty \|
\Delta^\rho_{\alpha_j} f - D^\rho f \|_p
&& \text{by Proposition~\ref{analysisleb}} \\
& \to 0 \qquad \text{as $j \to \infty$.}
\end{align*}
For $B_j$, we first note that if $P\eta=1$ a.e.\ then $P\eta_\rho =
1$ a.e. Thus Theorem~\ref{sample} implies in part (a) that $B_j \to
D^\rho f$ in $L^p$. In part (b), Theorem~\ref{sample} implies that
if the dilations expand exponentially, then $\frac{1}{J}
\sum_{j=1}^J B_j \to D^\rho f$ in $L^p$, since $\int_\Rd \eta_\rho
\, dx = \int_\Rd \eta \, dx = 1$. This proves the theorem.
\end{proof}

\subsection{\bf Proof of Theorem~\ref{sobolevrep2} --- synthesis onto
$W^{m,p}$} \ \label{sobolevrep2_proof}

In proving $S$ is surjective, we can assume the dilations expand
exponentially, like we did in the proof of
Theorem~\ref{lebesguerep2}.

First we prove surjectivity assuming $b=I$. Let $f \in W^{m,p}$
and $J \in \N$. Like in the proof of Theorem~\ref{lebesguerep2}
(but with $b=I$), we take $\phi=\charfn_{\cube}$ and define $c_J$
by \eqref{cJdef}. Then
\begin{equation} \label{srep2eq1}
\| c_J \|_{\ell^1(w^{m,p},\alpha)} \leq \| f \|_{W^{m,p}}
\end{equation}
by using Proposition~\ref{analysissob} for each $j=1,\ldots,J$. And
$S_{\psi,I} c_J \to f$ in $W^{m,p}$ as $J \to \infty$,  by
Theorem~\ref{scaleaveragesobolev}(b).

Hence the hypotheses of the open mapping theorem in
Appendix~\ref{banachapp} are satisfied with $A=1$, by
\eqref{srep2eq1}. Therefore $S_{\psi,I} : \ell^1(w^{m,p},\alpha) \to
W^{m,p}$ is open and surjective, and for each $f \in W^{m,p}$ and
$\e>0$ there exists $c \in \ell^1(w^{m,p},\alpha)$ with $S_{\psi,I}
c = f$ and $\| c \|_{\ell^1(w^{m,p},\alpha)} \leq \| f \|_{W^{m,p}}
+ \e$. This proves the theorem when $b=I$.

For the general case where $b \neq I$, we rescale like in the
proof of Theorem~\ref{sobolevrep1}: since $M_b \psi$ and $M_b
\eta$ satisfy formula \eqref{psidef} with ``$b=I$'', the case of
the theorem already proved tells us that for each $f \in W^{m,p}$
and $\e>0$, a sequence $c \in \ell^1(w^{m,p},\alpha)$ exists such
that $\| c \|_{\ell^1(w^{m,p},\alpha)} \leq \| f \|_{W^{m,p}} +
\e$ and $f = S_{M_b \psi,I} c$. The calculation \eqref{rescaling}
now implies $M_b^{-1} f = S_{\psi,b} c$, as desired.

\subsection{\bf Proof of Corollary~\ref{sobolevnorm3} --- analysis norm for $W^{m,p}$}\
\label{sobolevnorm3_proof}

By boundedness of the analysis operator in
Proposition~\ref{analysissob},
\[
\sup_j \sum_{|\rho|\leq m} |\alpha_j|^{|\rho|} \| \Delta^{\rho} T_j
f \|_{\ell^p} \leq \|P|\phi|\|_\infty \|f\|_{W^{m,p}}.
\]
On the other hand, choosing $\eta$ as in
Theorem~\ref{scaleaveragesobolev}(a), we see that $S_j T_j f \to f$
in $W^{m,p}$ as $j \to \infty$. It then follows from
Theorem~\ref{sobolevrep1} (bounded synthesis) that
\begin{align*}
\|f\|_{W^{m,p}}
& \leq \sup_j \| S_j T_j f \|_{W^{m,p}} \\
& \leq \|P|\eta|\|_{L^p(\cube)} \sup_j \sum_{|\rho| \leq m}
|\alpha_j|^{|\rho|} \| \Delta^\rho T_j f \| _{\ell^p} .
\end{align*}

\section*{\bf Acknowledgments} We thank Joaquim Bruna for showing
us his work \cite{B05} in preprint form, and we thank Maciej
Paluszy\'{n}ski and Guido Weiss for discussing the discrete
origins of the Hilbert transform with us.

\appendix

\section{\bf The open mapping theorem}
\label{banachapp}

The open mapping theorem in the following form is used to prove
surjectivity of the synthesis operator, at various points in the
paper.
\begin{theorem} \label{banach}
Let $X$ and $Y$ be Banach spaces, and suppose $S : X \to Y$ is
bounded and linear. Assume
\[
\overline{S(B_X(A))} \supset B_Y(1)
\]
for some $A>0$. That is, assume for each $y \in Y$ that a sequence
$\{ x_J \} \subset X$ exists with $Sx_J \to y$ as $J \to \infty$ and
$\| x_J \|_X \leq A \| y \|_Y$ for all $J$.

Then $S$ is an open mapping, and $S(X)=Y$. Indeed, given $y \in Y$
and $\e>0$ there exists $x \in X$ with $Sx=y$ and $\| x \|_X \leq A
\| y \|_Y + \e$.
\end{theorem}
For a proof, see \cite[Theorem~4.13]{rud} with $A=1/\delta$.

The hypothesis in Theorem~\ref{banach} can clearly be weakened, to
assume only for some \emph{dense} subset of $y$-values that a
sequence $\{ x_J \} \subset X$ exists with $Sx_J \to y$ as $J \to
\infty$ and $\| x_J \|_X \leq A \| y \|_Y$ for all $J$.

\section{\bf Discrete Hardy spaces}
\label{hardyapp}

In this appendix we study properties of the discrete Hardy space
$h^1$, which was defined in Section~\ref{hardyresults}. We will show
$h^1$ is independent of the cut-off function used in its definition,
and that it coincides (when $b=I$) with the discrete Hardy space
$H^1(\Zd)$ studied by Q. Y. Sun \cite{S93} and C. Eoff \cite{E95} in
dimension $1$, and later by S. Boza and M. J. Carro \cite{BC98,BC02}
in all dimensions.

We will need the following result on Riesz transforms of Schwartz
functions, which is a special case of \cite[Corollary~2.4]{BC02}.
Recall that $Z(x)=C_d x/|x|^{d+1}$ for $x \neq 0$, and $Z(0)=0$.
\begin{lemma} \label{rieszestimate}
If $\theta$ is a Schwartz function then
\[
R\theta(x) = \widehat{\theta}(0) Z(x) + O \left( \frac{1}{|x|^{d+1}}
\right) \qquad \text{for all $x\in \Rd \setminus \{ 0 \}$.}
\]
\end{lemma}
%

The next theorem replaces the kernel $z$ defining $h^1$ with a
discretization of the Riesz kernel, namely the sequence $z^b= \{
Z(bk) \}$ having $k$th term $z^b_k = Z(bk)$.
\begin{proposition}\label{hardyprop}
In the definition of the space $h^1$, if we replace the sequence
$z=\{z_k\}$ by $z^b=\{Z(bk)\}$ then we obtain the same space, with
an equivalent norm.
\end{proposition}
\begin{proof}[Proof of Proposition~\ref{hardyprop}]
Let $\nu$ be a cut-off function as in Section~\ref{hardyresults}.
Let $\mu$ be the Schwartz function with $\widehat{\mu}=\nu$. Note
$\widehat{\mu}(0)=\nu(0)=1$. Then
\[
\zeta(\xi) = -i \frac{\xi b^{-1}}{|\xi b^{-1}|} \nu(\xi) =
\widehat{Z}(\xi b^{-1}) \widehat{\mu}(\xi) = \widehat{K\mu}(\xi) ,
\qquad \xi \in \cube_0 ,
\]
where $K$ is the singular integral operator with kernel
\[
K(x) = |\det b| Z(bx) = |\det b| C_d \frac{bx}{|bx|^{d+1}} .
\]
Thus by definition of the sequence $z$ in
Section~\ref{hardyresults}, we have $z_k = \widehat{\zeta}(-k)=
K\mu(k)$.

Next, observe that for each $x\in \Rd$,
\begin{align*}
K\mu(x)
& = |\det b| C_d \, \text{p.v.} \! \int_{\Rd}\frac{by}{|by|^{d+1}} \mu(x-y) \, dy \\
& = C_d \, \text{p.v.} \! \int_{\Rd} \frac{y}{|y|^{d+1}} \mu(b^{-1}(bx-y)) \, dy \\
& = R\theta(bx),
\end{align*}
where $\theta(y)=\mu(b^{-1}y)$ is a Schwartz function with
$\widehat{\theta}(0) \neq 0$. In particular
\begin{align*}
z_k = K\mu(k) = R\theta(bk) & = \widehat{\theta}(0) Z(bk) +
O \left(\frac{1}{|bk|^{d+1}} \right) && \text{by Lemma~\ref{rieszestimate}} \\
& = \widehat{\theta}(0) z^b_k + O \left( \frac{1}{|k|^{d+1}} \right)
\end{align*}
for all $k\neq 0$. Hence $z-\widehat{\theta}(0) z^b \in \ell^1$, so
that $z$ and $z^b$ define identical $h^1$ spaces with equivalent
norms.
\end{proof}

The proposition and its proof yield a large class of kernel
sequences that generate $h^1$, for they show that if $\mu$ is any
Schwartz function with $\widehat{\mu}(0) \neq 0$, then we can
replace the sequence $\{ z_k \}$ by $\{ K \mu(k) \}$ in the
definition of $h^1$.

Now we can show independence of $h^1$ from the cut-off function.
\begin{corollary} \label{hardycor1}
The space $h^1$ does not depend on the cut-off function $\nu$ used
to define it, and different cut-off functions produce equivalent
norms.
\end{corollary}
\begin{proof}[Proof of Corollary~\ref{hardycor1}]
This follows from Proposition~\ref{hardyprop}, because $z^b$ does
not depend on $\nu$.

Alternatively, consider two different cut-off functions $\nu_1$ and
$\nu_2$, giving rise to periodic functions $\zeta_1$ and $\zeta_2$
as in Section~\ref{hardyresults}. Then $\zeta_1 - \zeta_2$ is smooth
and compactly supported in $\cube_0$ and hence has Fourier
coefficients in $\ell^1$. Therefore the kernel sequences associated
with $\nu_1$ and $\nu_2$ differ by only an $\ell^1$ sequence, and so
they define the same $h^1$ space, with comparable $h^1$ norms.
\end{proof}

\begin{corollary} \label{hardycor2}
If $b=I$ then $h^1=H^1(\Zd)$, with equivalent norms.
\end{corollary}
This corollary simply restates Proposition~\ref{hardyprop} with
$b=I$, because the discrete Hardy space $H^1(\Zd)$ is defined
(following \cite{BC02}) by the kernel sequence $z^I = \{ Z(k) \}$;
in other words,
\[
H^1(\Zd) = \{s \in \ell^1 : z^I * s \in \ell^1 \}
\]
with a norm $\| s \|_{H^1(\Zd)} = \| s \|_1 + \| z^I * s \|_1$.
Note that in dimension $1$, the sequence $z^I = \{ \charfn_{ \{ k
\neq 0 \} }/\pi k \}$ is called the \emph{Hilbert sequence} and
was considered by R. E. Edwards and G. I. Gaudry \cite{EG77}, who
proved boundedness of $s \mapsto z^I * s$ on $\ell^p(\Z)$, for
$1<p<\infty$.

S. Boza and M. J. Carro \cite{BC02} proved the space $H^1(\Zd)$
admits a characterization by maximal functions in the sense of
Fefferman--Stein \cite{FS72}, and an atomic decomposition in the
sense of Coifman--Weiss \cite{CW77}. The atomic decomposition in one
dimension was also stated in \cite{CW77}. It is an interesting
problem to investigate these characterizations for our space $h^1$
when $b$ is not the identity matrix.

\vspace{6pt} \noindent \emph{Remark on vanishing means in $h^1$.} If
$s \in h^1$ then $\sum_{k \in \Zd} s_k = 0$. \emph{Proof:} Writing
$\sigma(\xi) = \sum_{k \in \Zd} s_k e^{-2\pi i\xi k}$ we see that
$\sigma(\xi)\zeta(\xi)=\sum_{k \in \Zd} (s*z)_k e^{-2\pi i\xi k}$ in
$L^2(\cube_0)$. This last series is continuous because $s*z \in
\ell^1$, and $\sigma(\xi)$ is continuous too. But $\zeta(\xi)$ is
not continuous at the origin and so $\sigma(0)$ must equal zero, as
claimed.

Conversely if $b=I$ and $s \in \ell^1$ is finitely supported with
$\sum_{k \in \Zd} s_k = 0$, then $s \in h^1$; cf.\
\cite[Theorem~3.3]{BC02}. In other words, atoms belong to
$H^1(\Zd)$.

\vspace{6pt} We end this appendix with a question: is $h^1$
independent of the choice of ``translation'' matrix $b$? We suspect
not. Of course there is a trivial result: one can always replace $b$
by a multiple of $b$ without affecting the resulting space $h^1$.

\section{\bf Banach frames} \label{bframe}

This appendix explains how Banach frames arise from the analysis
norms earlier in the paper.

Let $Y$ be a Banach space, and let $Z$ be a Banach space whose
elements are complex sequences indexed by a countable set $I$. Let
$\{g_i\}_{i\in I}$ be a subset of $Y^*$, the dual space of $Y$, and
let $S_* : Z \to Y$ be a bounded linear operator. We say that
$(\{g_i\},S_*)$ is a {\em Banach frame} for $Y$ with respect to $Z$
if the following three conditions are satisfied:
\begin{itemize}
\item[(i)] $\{\langle f,g_i\rangle\} \in Z$, for all $f\in Y$,
\item[(ii)] $\|f\|_Y \approx \|\{\langle f,g_i\rangle\}\|_Z$, for all
  $f \in Y$,
\item[(iii)] $S_*(\{\langle f,g_i\rangle\}) = f$ for all $f\in Y$.
\end{itemize}
In other words, ``analyzing'' with the $\{ g_i \}$ maps $Y$ to $Z$
with comparable norms, and then ``synthesizing'' with $S_*$ recovers
the identity map on $Y$. The above definition is due to K.
Gr\"{o}chenig \cite{Gro91}; see the treatment in
\cite[\S17.3]{Ch03}.

The next result reformulates our $L^p$-analysis norm in
Corollary~\ref{lebesguenorm3} as a Banach frame result.
\begin{corollary}[Banach frame for $L^p$]
\label{banachlebesgue} Assume $1 \leq p < \infty$ and let $\phi$ and
$\psi$ satisfy the assumptions of Theorem~\ref{sample} with
$\gamma=1$ (scale-averaged approximation in $L^p$). Assume the
dilations expand exponentially.

Then $(\{ |\det b| \phi_{j,k}\},S_*)$ is a Banach frame for $L^p$
with respect to
\[
Z = \{ c \in \ell^\infty(\ell^p) : S_* c = \lim_{J \to \infty}
\frac{1}{J} \sum_{j=1}^J \sum_{k \in \Zd} c_{j,k} \psi_{j,k}\;
\text{exists in $L^p$} \}.
\]
\end{corollary}
One can similarly reformulate the Hardy and Sobolev space results
(Corollaries~\ref{hardynormb} and \ref{sobolevnorm3}), for isotropic
dilations $a_j=\alpha_j I$.
\begin{proof}[Proof of Corollary~\ref{banachlebesgue}]
Write $c_j = \{ c_{j,k} \}_{k \in \Zd}$, so that $S_j c_j = \sum_{k
\in \Zd} c_{j,k} \psi_{j,k}$. Define $Bc=\{ S_j c_j \}_{j>0}$. Then
$B$ is a bounded linear operator from $\ell^\infty(\ell^p)$ into
$\ell^\infty(L^p)$ (by the proof of Theorem~\ref{lebesguerep1}).
Therefore Lemma~\ref{seqspace} below tells us $Z$ is a closed
subspace of $\ell^\infty(\ell^p)$, and hence is a Banach space
itself. Moreover, $S_*$ is bounded from $Z$ into $L^p$.
Corollary~\ref{lebesguenorm3} and Theorem~\ref{sample} now show that
$( \{ |\det b| \phi_{j,k} \}, S_*)$ is a Banach frame for $L^p$ with
respect to $Z$, noting in particular by Theorem~\ref{sample} that
\[
S_*(Tf) = \lim_{J \to \infty} \frac{1}{J} \sum_{j=1}^J S_j T_j f = f
\qquad \text{in $L^p$, for all $f \in L^p$.}
\]
\end{proof}

The final lemma states that the preimage of the space of
Ces\`{a}ro-convergent sequences in a Banach space forms a closed
subspace. Let $Y$ be a Banach space, and write $\ell^\infty(Y)$ for
the Banach space of all sequences $y = \{y_j\}_{j>0}, y_j \in Y$,
such that $\| y \|_{\ell^\infty(Y)} = \sup_{j>0} \| y_j \|_Y <
\infty$.
\begin{lemma} \label{seqspace}
Let $X$ and $Y$ be Banach spaces and $B: X \to \ell^\infty(Y)$ be a
bounded linear operator. Then the subspace
\[
\{ x \in X : \text{the limit $\lim_{J \to \infty} \frac{1}{J}
\sum_{j=1}^J (Bx)_j$ exists in $Y$} \}
\]
is closed in $X$.
\end{lemma}
We omit the proof.

\end{document}